\def\obs#1{{\bf (*** #1 ***)} }

\def\obs#1{}     % Remova esta linha para rodar a versao 1
\documentclass[twoside,letterpaper,draft,11pt]{amsart}
\usepackage{amsmath,amsthm,latexsym,xspace,amscd,amssymb,mathtools,color,enumerate}

\DeclareUnicodeCharacter{00B4}{}

\usepackage[utf8]{inputenc}

\usepackage[mathscr]{eucal}
\usepackage{amssymb}
\usepackage[all]{xy}
\makeatletter\renewcommand\theenumi{\@roman\c@enumi}\makeatother

\newtheorem{teo1}{Theorem}[section]
\newtheorem{def1}[teo1]{Definition}
\newtheorem{lem1}[teo1]{Lemma}
\newtheorem{cor1}[teo1]{Corollary}
\newtheorem{prop1}[teo1]{Proposition}

\newtheorem{exe}[teo1]{Example}
\newtheorem{remark}[teo1]{Remark}

\newcommand{\m}{{}^{-1}}

\newcommand{\U}{\mathcal U}

\newcommand{\id}{{\rm id}}
\newcommand{\Z}{{\mathbb Z}}

\newcommand{\N}{{\mathbb N}}

\newcommand{\A}{{A}}

\newcommand{\de}{\delta}

\newcommand{\af}{\alpha}
\newcommand{\bt}{\beta}

\newcommand{\gm}{\gamma}

\newcommand{\ta}{\tau}

\DeclareMathOperator{\G}{G}

\def\ndv{\ {\mid \kern -0.7 em {\scriptstyle \not}} \ \ }

\def\nd{\ {\mid \kern -0.4 em {\scriptstyle \not}} \ \ }

\begin{document}

\thispagestyle{empty}

\title[On partial skew groupoid rings]{Ring theoretic properties of partial skew groupoid rings with applications to Leavitt path algebras}

\author[D. Bagio]{Dirceu Bagio}
\address{ Departamento de Matem\'atica, Universidade Federal de Santa Maria, 97105-900\\
	Santa Maria-RS, Brasil}
\email{bagio@smail.ufsm.br}

\author[V. Mar\'in ]{V\'ictor Mar\'in}
\address{Departamento de Matem\'{a}ticas y Estad\'{i}stica, Universidad del Tolima, Santa Helena\\
	Ibagu\'{e}, Colombia} \email{vemarinc@ut.edu.co}

\author[H. Pinedo]{H\'ector Pinedo}
\address{Escuela de Matematicas, Universidad Industrial de Santander, Cra. 27 Calle 9 ´ UIS
	Edificio 45, Bucaramanga, Colombia}
\email{hpinedot@uis.edu.co}


\date{\today}
\begin{abstract}
	Let $\alpha=(\A_g,\af_g)_{g\in \G}$ be a group-type partial action of a  connected grou\-poid $\G$ on a ring $A=\bigoplus_{z\in \G_0}A_z$ and $B:=A\star_{\af}\G$ the corresponding partial skew groupoid ring. In the first part of this paper we investigate the relation of several ring theoretic properties between $A$ and $B$. For the second part, using that every Leavitt path algebra is isomorphic to a partial skew groupoid ring obtained from a partial groupoid action $\lambda$, we characterize when $\lambda$ is group-type. In such a case, we obtain ring theoretic properties of Leavitt path algebras from the results on general partial skew groupoid rings. Several examples that illustrate the results on Leavitt path algebras are presented. 
\end{abstract}

\maketitle

\setcounter{tocdepth}{1}

\section{Introduction}
\label{intro}
Let $\G$ be a groupoid and $A$ an associative and unital ring. For a partial action $\af$ of $\G$ on $A$ corresponds the partial skew groupoid ring $B:=A\star_\alpha \G$ (cf. Definition \ref{def-partial-skew}). The partial skew groupoid ring is a generalization of both partial skew group ring and groupoid algebra. Precisely, if $\G$ is a group then $B$ is the partial skew group ring while if $\alpha$ is the trivial partial action $\alpha_g=\id_A$, for all $g\in G$, then $B$ is the usual groupoid algebra. \smallbreak

Partial skew groupoid rings were firstly considered in \cite{BFP} for ordered grou\-poids and, after that, in \cite{BP} for general groupoids. Basic properties of $B$ were given in \cite{BFP,BP}. For instance, if $\af=(\A_g,\af_g)_{g\in \G}$ then $B$ is unital if and only if $A_z$ is unital, for all $z\in \G_0$ and $\G_0$ is  finite, where  $\G_0$ denotes the set of objects of $\G.$  Also, if $A_g$ is a unital ring for each $g\in \G$, then $B$ is an associative ring.   
Partial skew groupoid rings appear naturally in the partial Galois theory of groupoids \cite{BP} and in the context of groupoid graded rings \cite{NyOP2020}. Also, it was proved in \cite{GY} that every  Leavitt path algebra is a partial skew groupoid ring.  \smallbreak

A groupoid $\G$ is a disjoint union of its connected components which are full connected subgroupoids of $\G$. We use this decomposition of $\G$ to reduce  partial groupoid actions to the connected case; see Remark 3.3 of \cite{BPP}. The structure of a connected groupoid is well-known. Any connected groupoid $\G$ is isomorphic to $\G_{0}^2\times \G(x)$, being $\G(x)$  the isotropy group of an object $x\in \G_0$. If $A=\bigoplus_{y\in \G_0}A_y$ and $\alpha$ is a group-type partial action (see Definition \ref{def-group-type}) then the factorization of $\G$ induces a factorization of $B$. In fact, it was proved in \cite{BPP} that there are a global action $\beta$ of $\G_0^2$ on $A$ and  a partial action $\gamma$ of $\G(x)$ on $A\ast_{\beta}\G_0^2$  such that $B\simeq \big(A\ast_{\beta}\G_0^2\big)\ast_{\gamma}\G(x)$. In other words, if $\alpha$ is group-type then the partial skew groupoid ring $B$ is a partial skew group ring.   \smallbreak

In the first part of this paper we investigate ring theoretic properties of $B$ for the case that it admits the factorization $B\simeq \big(A\ast_{\beta}\G_0^2\big)\ast_{\gamma}\G(x)$. 
The properties of $B$ that are studied are the following: noetherianity, von Neumann regularity, the Jacobson radical, perfect, semiprimary and Krull dimension. One of the strategies used to prove the results is to apply the ring theoretic properties proved in \cite{BC} to $B$ since it is a skew partial group ring. \smallbreak

In the second part, we work with Leavitt path algebras. Precisely, let $L_{\Bbbk}(E)$ be the Leavitt path algebra associated to a directed graph $E$ over a field $\Bbbk.$ By \cite{GY}, we have that $L_{\Bbbk}(E)\simeq D(X)\star_{\lambda}\G(E)$, where $D(X)$ is $\Bbbk$ algebra and $\lambda$ is a partial action of the free path groupoid $\G(E)$ on $D(X)$. A general characterization for $\lambda$ to be group-type is given in Theorem \ref{gtype}. A refinement of this result is given in Proposition \ref{daniel} for the case that $E$ has a sink, $E^0$ is finite and $|E^{0}|\geq 2$. Assuming that $\lambda$ is group-type, we present ring theoretic properties of $L_{\Bbbk}(E)$. Examples to illustrate the results on Leavitt path algebras are also presented. \smallbreak

The paper is organized as follows. The background on groupoids and partial groupoid actions will be given in Section 2. The results on ring theoretic properties of the partial skew groupoid ring are presented in Section 3. Parti\-cu\-larly, for noetherianity and regularity properties it is convenient to establish first some general results of groupoid graded rings. For this, we introduce the notion of von Neumann groupoid graded ring and, in  Proposition \ref{gvng}, we relate this concept with von Neumann regularity. 
The results related to Leavitt path algebras are given in Section 4.

\subsection*{Conventions}\label{subsec:conv}
Throughout this work,  by ring we mean an associative ring with identity. A ring with identity will be called unital. The Jacobson radical of a ring $A$ will be denoted by $J(A)$. The group of invertible elements of a unital ring $A$ is denoted  by $\mathcal{U}(A)$ and  the characteristic of a field $\Bbbk$ by $\operatorname{ch}(\Bbbk)$. 

\section{Preliminaries} 

\subsection{ On connected groupoids} A {\it groupoid} $\G$ is a small category in which every morphism is an isomorphism. The set of  objects of $\G$ is  denoted by $\G_0$. The {\it source} and {\it the target maps} $s,t:\G\to \G_0$ are defined by $s(g)=x$ and $t(g)=y$, for any morphism $g:x\to y$ of $\G$. The set $\G(x,y)$ contains the morphisms $g$ of $\G$ for which $s(g)=x$ and $t(g)=y$. Every object $x$ of $\G$ is identified with the identity morphism, that is, $x=\id_x$. Thus $\G_0\subseteq \G,$ $s(g)=g^{-1}g$ and $t(g)=gg\m$, for all $g\in \G$. For each $x\in \G_0$, the {\it isotropy group} associated to $x$ is  $\G(x):=\G(x,x)$ and we will fix $\G(x,\,):=\{g\in \G\,:\,s(g)=x\}$ and  $\G(\,,x):=\{g\in \G\,:\,t(g)=x\}$.\smallbreak

We denote the composition of morphisms of a groupoid by concatenation. Notice that, for $g,h\in \G$, there exists $gh$ if and only if $t(h)=s(g)$. In this case, we say that $g$ and $h$ are composable. We  denote by $\G_s\times _t\!\G$ the subset of $\G\times \G$ of composable pairs, that is, $\G_s\times _t\!\G$ is the pullback of $s$ and $t$ in the category of groupoids. Observe that $s(gh)=s(h)$ and $t(gh)=t(g)$, for all $(g,h)\in \G_s\times _t\!\G$.     \smallbreak

A groupoid $\G$ is called {\it connected} if $\G(x,y)\neq \emptyset$, for any $x,y\in \G_0$.  For an  arbitrary groupoid  $\G$  the equivalence relation on $\G_0$, given by $x\sim y$  if and only if $\G(x,y)\neq \emptyset$,  
induces a decomposition of $\G$ as a disjoint union of connected subgroupoids. Explicitly, for each equivalence class $X\subset\G_0$ corresponds the full subgroupoid $\G_X$ of $G$ whose set of objects is $X$. Hence $\G$ is the disjoint union of the subgroupoids $\G_X$, with $X\in \G_0/\!\!\sim$. \smallbreak

In order to present the structure of a connected groupoid, we recall that the  {\it coarse groupoid} $Y^2$ (or Brandt groupoid) associated to any non-empty set $Y$ is the groupoid whose objects are the set $Y$ and satisfies $s(y,y')=y$, $t(y,y')=y'$ and $(y',y'')(y,y')=(y,y'')$. The next result is well known and its proof will be omitted; see, for instance, Proposition 2.1 of \cite{BPP}.

\begin{prop1}\label{group:connec}
	Let $\G$ be a connected groupoid. Then $\G\simeq \G_0^2\times \G(x)$ as groupoids.
\end{prop1}

\subsection{Partial groupoid action}
Throughout this subsection, $\G$ denotes a groupoid and $A$ denotes a ring. We start by recalling the definition of partial action given in \cite{BP}.
\begin{def1}\label{def-groupoid-action}
	{\rm  A \emph{partial action} of $\G$ on $\A$ is a set of pairs $\af=(\A_g,\af_g)_{g\in \G}$  that satisfies:
		\begin{enumerate}[\rm (i)] 
			\item  $\A_g$ is an ideal of $\A_{t(g)}$, $\A_{t(g)}$ is an ideal of $\A$ and $\af_g:\A_{g\m}\to \A_g$ is a ring isomorphism, for all $g\in \G$,\vspace{.05cm}
			\item  $\alpha_x=\id_{A_x}$, for all $x\in \G_0$,\vspace{.05cm}
			\item $\af_h^{-1}(\A_{g^{-1}}\cap\A_h)\subset \A_{(gh)^{-1}}$, for all $(g,h)\in \G_s\times _t\!\G$,\vspace{.05cm}
			\item $\af_{gh}(a)=\af_g\af_h(a), \, \text{ for all }\, a\in \af_h^{-1}(\A_{g^{-1}}\cap\A_h)$, $(g,h)\in \G_s\times _t\!\G$.
	\end{enumerate}}
\end{def1}

\smallbreak

We also recall other definitions that will be used later.

\begin{def1}\label{def-specific}
	{\rm A partial action $\af=(\A_g,\af_g)_{g\in \G}$ of a groupoid $\G$ on a ring $A$ is called
		\begin{enumerate}[\rm (i)] 
			\item  {\it global}, if $\af_g\af_h=\af_{gh}$, for all  $(g,h)\in\G_s\times _t\!\G$,
			\item  {\it unital} if each $A_g$ is a unital ring, that is, there exists a central idempotent element $1_g$ in $A$ such that $A_g=A1_g$, for all $g\in \G$,
			\item  {\it  finite-type}, if for any  
			$z\in \G_0$ there are $g_1,\ldots, g_n\in \G(\,,z)$ such that 
			$A_{t(g)}=\sum_{i=1}^n A_{gg_i},$ for any $g\in \G(z,\,).$
	\end{enumerate}} 
\end{def1}

For the convenience of the reader, we recall Lemma 1.1 of \cite{BP} which gives some useful properties of partial groupoid actions.

\begin{lem1}\label{lem:BP} Let $\af=(\A_g,\af_g)_{g\in \G}$ be a partial action of a groupoid $\G$ on a ring $A$. The following statements are true:
	\begin{enumerate}[\rm (i)]
		\item  $\af_{g\m}=\af\m_g$, for all $g\in \G$,\vspace{.05cm}
		\item $\af_g(\A_{g\m}\cap \A_h)=\A_{g}\cap \A_{gh}$, for all $(g,h)\in \G_s\times _t\!\G$,\vspace{.05cm}
		\item  $\alpha$ is global if and only if $A_g=A_{t(g)}$, for all $g\in \G$.
	\end{enumerate}
\end{lem1}

\begin{remark} \label{obs-pag}{\rm  Let $\af=(\A_g,\af_g)_{g\in \G}$ be a partial action of a groupoid $\G$ on a ring $A$. Note  that $\af$ induces a partial action $\af_{\G(x)}=(\A_h,\af_h)_{h\in\G(x)}$ of the isotropy group $\G(x)$ on the ring $\A_x$, for each $x\in\G_0$.  Moreover, using the decomposition of $\G=\dot\cup_{X\in \G_0/\!\sim}\G_X$ in connected components, it follows by Proposition 2.3 of \cite{BPi} that  partial actions of $\G$ on $A$ are uniquely determined by partial actions of the connected groupoids  $\G_X$, $X\in \G_0/\!\!\sim$, on $A$. }
\end{remark}

Let $\G$  be a connected groupoid and $x\in \G_0$. Consider the following equivalence relation on $\G(x,\,)$:
$$g\sim_x l\,\,\ \text{ if and only if }\,\, t(g)=t(l).$$
A transversal $\tau(x)=\{\tau_{y}\mid y\in \G_0\}$ for $\sim_x$ such that $\tau_x=x$ will be called a \emph{transversal for $x$}. Observe that $\tau_y\in\G(x,y)$, for all $x\neq y\in\G_0$, and $\ta_x=x\in \G(x)$. The next definition was given in subsection 3.2 of \cite{BPP}.
\begin{def1}\label{def-group-type}
	{\rm A partial action $\af=(A_g,\af_g)_{g\in \G}$ of a connected groupoid $\G$ on a ring $A$ will be called {\it  group-type} if there exist $x\in \G_0$ and a transversal $\tau(x)=\{\tau_{y}\mid y\in \G_0\}$ for $x$ such that
		\begin{align}\label{cond1} 
			A_{\tau\m_y}=A_x \ \ \text{and} \ \ A_{\tau_y}=A_y, \ \ \text{ for all } \ y\in\G_0.
	\end{align}}
\end{def1}

\begin{remark}\label{obs-not-depend}
	{\rm	\begin{enumerate}[\rm (i)]
			\item By Lemma \ref{lem:BP} (iii), any global groupoid action is group-type. Also, the notion of group-type partial action does not depend on the choice of object $x$, cf. Remark 3.4 of \cite{BPP}. In fact, if $\ta(x)$ is a transversal for $x$ and $y\in \G_{0}$ then $\gamma(y):=\{\tau_z\tau\m_y\mid z\in \G_0\}$ is a transversal for $y$. \smallbreak
			
			\item The name group-type is suggested by Corollary 3.6 of  \cite{BaPaPi} in which it is observed that this kind of partial actions  can be obtained  from partial group actions. 
	\end{enumerate}}
\end{remark}

Assume that $\af=(A_g,\af_g)_{g\in \G}$ is a group-type partial action of a connected groupoid $\G$ on a ring $A$. From \eqref{cond1} follows  $A_{g\m}=A_{g\m}\cap A_{s(g)}= A_{g\m}\cap A_{\ta_{s(g)}}$ which implies that 
$A_g=\af_g(A_{g\m})= \af_g(A_{g\m}\cap A_{\ta_{s(g)}})=A_g\cap A_{g\ta_{s(g)}}$, for any $g\in \G$. Hence
\begin{align}\label{eq1}
	A_g\subseteq A_{g\tau_{s(g)}},\,\,\, g\in \G.
\end{align}

\begin{remark} \label{ob-invar-ideal}
	{\rm Let $\af$ be a partial action of a groupoid $\G$ on a ring $A$. Consider an  $\alpha$-invariant ideal $I$ of $A$, that is, $\af_g(I\cap A_{g\m})\subseteq I$, for all $g\in \G$. In this case the family $\af|_I=(I_g, \af_g)_{g\in \G}$, where $I_g:=I\cap A_{g}$, is a partial action of $\G$ on $I$. It is clear that if $\af$ is group-type then $\af|_I$ is group-type.}
\end{remark}

%\begin{exe}\label{ex:groupoid-12}{\rm
%		Let $\G=\{g,g^2,h,h^2,l,l\m,m,m\m,n,n\m\}$ be the groupoid with set of objects $\G_0=\{x,y\}$. The composition in $\G$ is given by
%		\begin{align*}
%			g&\in \G(x),\quad h\in \G(y),& &\,\,\,g^3=x,& &\,h^3=y,& \\
%			&l,m,n\in\G(x,y),& lg&=m=hl,& mg&=n=hm.&
%		\end{align*} 
%		The groupoid $\G$ is illustrates in the following diagram 
%		\[\xymatrix{& x\ar[d]^{n}\ar[r]^{g}  &x\ar[d]^{m} \ar[r]^{g} &x\ar[d]^{l}\\
%			& y & y\ar[l]^{h}&\ar[l]^{h} y}\\ \]
%		Now we define a partial action of $\G$ on $A=$  {\color{red} TERMINAR.}	}
%
%\end{exe}

\subsection { The partial skew groupoid ring}
Throughout this subsection, $\af=(A_g,\af_g)_{g\in \G}$ denotes a unital partial action of a groupoid $\G$ on a ring $A$. We will assume that $A_g=A1_g$, where $1_g$ is a central idempotent of $A$, for all $g\in\G$.  
\begin{def1}\label{def-partial-skew}
	{\rm The {\it partial skew groupoid ring} $A\star_\af\G$ associated to $\alpha$ is the set of  formal sums $\sum_{g\in \G}a_g\delta_g$, where $a_g\in A_g$, with the usual addition and multiplication induced linearly by the following rule
		\begin{equation}\label{prod}
			(a_g\delta_g)(a_h\delta_h)=
			\begin{cases}
				a_g\alpha_g(a_h1_{g^{-1}})\delta_{gh}, &\text{if $(g,h)\in \G_s\times _t\!\G$,}\\
				0, &\text{otherwise},
		\end{cases}\end{equation}
		for all $g, h\in \G$, $a_g\in A_g$ and $a_h\in A_h$.}
\end{def1}

The partial skew groupoid ring $A\star_\af\G$ is an associative ring. It was proved in Section 3 of \cite{BFP} that if $\G_0$ is finite then $A\star_\af\G$ is  unital with identity  $1_{A\star_\af\G}=\sum_{y\in\G_0}1_y\delta_y$.  \smallbreak
%Notice that $A\star_\af\G$ is $\G$-graded by putting $(A\star_\af\G)_g=A_g\de_g$, for all $g\in \G$.

On the other hand it was shown in \cite{BPP} that, if $\G$ is connected then the factorization of $\G$ given by Proposition \ref{group:connec} induces, under suitable conditions, a factorization of the partial skew groupoid ring. In order to review this construction, assume  that $\G$ is connected, $\G_0$ is finite, $\af$ is group-type and $A=\bigoplus_{z\in \G_0}A_{z}$. We also fix $x\in \G_0$ and $\tau(x)=\{\tau_{y}\mid y\in \G_0\}$ a transversal for $x$ such that \eqref{cond1} is satisfied. \smallbreak

For each $u=(y,z)\in \G_0^2$, we consider 
\begin{align}\label{action-beta}
	&B_u=A_{t(u)}=A_z,& &\bt_u=\af_{\ta_{t(u)}}\af_{\ta\m_{s(u)}}:A_y\to A_z.& 
\end{align}
By Lemma 4.1 of \cite{BPP}, the family $\bt=(B_u,\bt_u)_{u\in\G_0^2}$ is a global action of $\G_0^2$ on $A$. Hence we may consider the skew groupoid ring $C:=A\star_ {\bt}\G_0^2$. There is a partial action $\gamma$ of $\G(x)$ on $C$ defined in the following way. Given $z\in \G_0$ and $h\in \G(x)$, set $C_{z,h}:=\af_{\ta_z}(A_h)$ and 
\begin{align}\label{Ch}
	C_h:=\bigoplus_{u\in\G_0^2}C_{t(u),h}\delta_u=\bigoplus_{u\in\G_0^2}\af_{\ta_{t(u)}}(A_h)\delta_u.
\end{align}
Also, $\gamma_{z,h}:C_{z,h\m}\to C_{z,h}$ is defined by $\af_{\ta_z}(a)\mapsto \af_{\ta_z}(\af_h(a))$, for all $a\in A_{h\m}$. These maps induce the following ring isomorphism
\begin{align}\label{gamma-h}
	\gamma_h:C_{h\m}\to C_h,\quad  \gamma_h(\af_{\ta_{t(u)}}(a)\delta_u)= \gamma_{t(u),h}(a)\delta_u, \,\,\, a\in A_{h\m},\,\, u\in \G_0^2.
\end{align}
From Lemmas 4.2 and 4.3 of \cite{BPP} follow that $\gamma=(C_h, \gamma_h)_{h\in \G(x)}$ is a unital partial action of $\G(x)$ on $C$. Moreover $C_h=C1'_h$, where $1'_h=\sum_{z\in \G_0}\af_{\tau_z}(1_h)\delta_{(z,z)}$, for all $h\in \G(x)$. In order to enunciate Theorem 4.4 of \cite{BPP}, which give us the factorization of $A\star_\af\G$, we consider the  notation: $g_x=\ta\m_{t(g)}g\ta_{s(g)}\in \G(x)$, for all $g\in \G$.

\begin{teo1}\label{teo-decomp} Suppose that $\G$ is connected, $\G_0$ is finite, $\af$ is group-type and $A=\bigoplus_{z\in \G_0}A_{z}$. Then
	the map $\psi:A\star_\af\G\to (A\star_{\bt}\G_0^2)\star_ \gamma\G(x)$  given by
	\begin{equation}\label{varphi} \psi(a\delta_g)=a\delta_{(s(g),t(g))}\delta_{g_x},\quad  a_g\in A_g,\,\, g\in \G, \end{equation}  is a ring isomorphism. 
\end{teo1}

The factorization of $A\star_\af\G$ given above will be useful in the rest of the paper. We end the background with the following.

\begin{lem1}\label{latoga} Suppose that $\af$ is group-type. If $\alpha$ is finite-type, then $\gamma$ defined in \eqref{gamma-h} is  finite-type.
\end{lem1}	
\begin{proof}
	Since $\alpha$ is finite-type, there are $g_1, \dots, g_n\in \G(\,,x)$ such that 
	$A_x=\sum_{i=1}^n A_{gg_i}$. For each $h\in \G(x)$ we take $\{hg_1\tau_{s(gg_1)}, \ldots, hg_n\tau_{s(gg_n)}\}\subseteq \G(x)$. Then 
	\begin{align*}
		A\star_{\bt}\G_0^2&=\bigoplus_{u\in \G_0^2}A_{t(u)}\delta_u= \bigoplus_{u\in \G_0^2}\alpha_{\tau_{t(u)}}(A_x)\delta_u=\bigoplus_{u\in \G_0^2}\sum_{i=1}^n\alpha_{\tau_{t(u)}}(A_{hg_i})\delta_u\\
		&=\sum_{i=1}^n\bigoplus_{u\in \G_0^2}\alpha_{\tau_{t(u)}}(A_{hg_i})\delta_u\stackrel{\eqref{eq1}}\subseteq\sum_{i=1}^n\bigoplus_{u\in \G_0^2}\alpha_{\tau_{t(u)}}(A_{hg_i\tau_{s(hg_i)}})\delta_u\\
		&\stackrel{\eqref{Ch}}=\sum_{i=1}^nC_{hg_i\tau_{s(hg_i)}}\subseteq A\star_{\bt}\G_0^2.
	\end{align*}
	Hence $A\star_{\bt}\G_0^2=\sum\limits_{i=1}^nC_{gg_i\tau_{s(gg_i)}}$ and we conclude that $\gamma$ is finite-type.
\end{proof}	

%
%\begin{remark}\label{obs-finite-type}{\rm
%If $\G$ is finite then $\alpha$ is finite-type. In fact, as $\alpha$ is unital, Theorem 2.1 of \cite{BP} implies that $\alpha$ is globalizable. Hence there exists a global action $\tilde{\alpha}=(\tilde{A}_g,\tilde{\alpha}_g)_{g\in \G}$ of $\G$ on a ring $\tilde{A}$ such that $\tilde{A}_z=\sum_{g\in \G(\,,z)}\tilde{\alpha}_g(A_z)$. From Proposition 6.4 (ii) of \cite{NOP} follows that $\af$ is finite-type.}
%\end{remark}

\section{Ring theoretic properties of $A\star_\af\G$}\label{sec-3}

In all what follows in this section, $\G$ is a connected groupoid such that $\G_0$ is finite and $\af=(A_g,\af_g)_{g\in \G}$ is a group-type unital partial action of $\G$ on  $A=\bigoplus_{z\in \G_0} A_z$. We will assume that there is $x\in \G_0$ and $\tau(x)=\{\tau_{y}\,:\,y\in\G_0\}$ a transversal for $x$ that  satisfies \eqref{cond1} and $A_g=A1_g$, where $1_g$ is a central idempotent of $A$, for all $g\in\G$. We also assume that $\beta$ and $\gamma$ are as in the previous section. \smallbreak

The aim in this section is to establish relations between ring theoretic properties of $A$ and $A\star_\af\G$.  First, observe that
\begin{align}\label{inclusion-map}
	\varphi: A\to A\star_\alpha \G\,\text{ given by }\,  \varphi(a)=\sum_{z\in \G_0}(a1_z)\delta_z,\quad a\in A,	
\end{align}
is an injective ring homomorphism and whence we identify $A$ as a unital  subring of $A\star_\alpha \G$. Moreover  $\varphi$ induces the following ring isomorphisms 
$A\simeq  \oplus_{z\in \G_0}A_z\delta_z\simeq \oplus_{w\in (\G^2_0)_0}A_w\delta_{w}$. Consequently
\begin{equation} \label{isosos}A\simeq (A\star_\af\G)_0 \simeq  (A\star_{\bt}\G_0^2)_0.
\end{equation}

\subsection{Noetherianity} 
The relation of the artinian property between $A$ and $A\star_\af\G$ was given in Theorem 1.3 of \cite{NOP}, for an arbitrary groupoid $\G$. In our context,  this result was refined in Theorem 5.11 of \cite{BPP}.  
Here we explore the relation of the noetherian property between $A$ and $A\star_\af\G$. \smallbreak

We recall some notions and facts on group and groupoid graded rings.  We say that  that a ring $R$ is called {\it $\G$-graded} if there is a set $\{ R_g \}_{g \in \G}$ of additive subgroups of
$R$ such that $R = \bigoplus_{g \in \G} R_g$, $R_gR_h\subset R_{gh}$ if $(g,h)\in \G_s\times _t\!\G$ and $R_g R_h = \{ 0 \}$ otherwise.
A $\G$-graded ring $R=\bigoplus_{g \in \G} R_g$ that satisfies $R_g R_h = R_{gh}$ for all $(g,h)\in \G_s\times _t\!\G$ is said {\it strongly graded}.  Then  by setting  $(A\star_\af\G)_g=A_g\delta_g$ for any $ g\in \G,$ follows from \eqref{prod} that $A\star_\af\G$ is a $\G$-graded ring. Moreover,  by Proposition 3.4 of \cite{BFP},  $A\star_\af\G$ is strongly graded if and only if $\alpha$ is global.  If $\G=G$ is a group we say that  $G$-graded  ring $R = \oplus_{g \in G} R_g$ is  {\it epsilon-strongly graded} if 
$R_g R_h = R_g R_{g^{-1}} R_{gh} = R_{gh} R_{h^{-1}} R_h,\quad \text{for all }g,h\in G.$ It is not difficult to see that unital partial crossed products are epsilon-strongly graded, see Theorem 35 of \cite{NyOP2018}  for details.

\begin{lem1} \label{noet} Let $\G$ be a finite groupoid, $R=\bigoplus_{g\in \G}R_g$ be a strongly $\G$-graded unital ring and $R_0=\bigoplus_{z\in \G_0}R_z$. If $R_0$ is left (right) noetherian, then $R$ is left (right) noetherian.
\end{lem1}

\begin{proof} 
	Since  $R_g R_{g\m}=R_{t(g)}$, for all $g\in \G$,  we obtain a $(R_g,R_{g\m})$-bimodule epimorphism  $\mu: R_g \otimes_{R_{s(g)}}R_{g\m}\to R_{t(g)}.$ Similarly, there exists a $(R_{g\m},R_{g})$-bimodule epimorphism $\tau: R_{g\m} \otimes_{R_{t(g)}}R_{g}\to R_{s(g)}$. Then, the Morita context  $(R_{t(g)}, R_{s(g)}, R_g, R_{g\m}, \mu, \tau)$ is strict. Thus ${R_g}$ is a left projective and finitely generated $R_{t(g)}$-module  and ${R_g}$ is a right projective and finitely generated $R_{s(g)}$-module. Consequently, ${R_g}$ is a left projective and finitely generated  right $R_0$-module . Using that $\G$ is finite we conclude that $R$ is a left and right  finitely generated $R_0$-module. The fact that $R_0$ is left (right) noetherian implies that  ${R}$ is a left and right noetherian $R_0$-module. Finally, from $R_0\subseteq R$ we conclude that $R$ is a left and right noetherian $R$-module  and the result follows. \end{proof}

In order to prove the main result of this subsection, we recall the following. A group $G$ is called a  {\it polycyclic-by-finite} if there exists a  subnormal series  $$\{1\}=G_0\unlhd G_1\unlhd \ldots\unlhd G_n\unlhd G_{n+1}=G$$ such that $G/G_n$ is finite and  $G_{i+1}/G_i$ is cyclic, for all $0\leq i\leq n-1$.

\begin{teo1}\label{noete} The following assertions hold.
	\begin{enumerate}[\rm (i)] 
		\item If  $A\star_\af\G$ is left  noetherian then $A$ and $A\star_{\bt}\G_0^2$ are left noetherian.\smallbreak
		
		\item If $A$ is left noetherian and $\G(x)$ is a polycyclic-by-finite group then   $A\star_\af\G$ is left noetherian.
	\end{enumerate}
\end{teo1}
\begin{proof}
	\noindent (i) From Theorem \ref{teo-decomp} and Proposition 3.4 of \cite{NOP} follow that $A\star_{\bt}\G_0^2$ is left  noetherian. Applying again Proposition 3.4 of \cite{NOP} for $A\star_{\bt}\G_0^2$ we conclude that $A$ is left noetherian. \smallbreak
	
	\noindent (ii) By Lemma \ref{noet} the ring  $C:=A\star_{\bt}\G_0^2$ is left noetherian. Since $C\star_ \gamma\G(x)$ is  epsilon strongly graded Theorem 3.7 of \cite{L} implies that it is left noetherian and whence the result follows from Theorem \ref{teo-decomp}.
\end{proof}

\subsection{Von Neumann regularity} Recall that a unital associative ring $R$ is von Neumann regular if $a\in aRa,$ for all $a\in R$.
There are several equivalent conditions to the notion of von Neumann regularity. Here, we will use the following equivalences: $R$ is von Neumann regular $\Leftrightarrow$ if every principal left ideal is a direct summand of the left $R$-module $R$ $\Leftrightarrow$ if every finitely generated left ideal is a direct summand of the left $R$-module $R$. \smallbreak

The next definition is inspired in the group case which was considered in Section 2.2 of \cite{La2}. A $\G$-graded ring  $R=\bigoplus_{g\in \G} R_g $ will be called a {\it graded von Neumann regular } if and only if $a\in aRa,$ for each $a\in R_g$ and $g\in \G.$ Clearly, a $\G$-graded ring that is von Neumann regular is graded von Neumann regular. The converse is not true even in the group case; see Example 2.4 of \cite{La2}.\smallbreak

\begin{remark}\label{obs-g-grade}{\rm 
		Let $B:=A\star_\af\G$ and $C:=A\star_{\beta}\G_{0}^2$. By Theorem \ref{teo-decomp}, we have $B\simeq C\star_{\gamma}\G(x)$. Notice that if $B$ is $\G(x)$-graded von Neumann regular then $B$ is $\G$-graded von Neumann regular. In fact, let $g\in \G$ and $a_g\delta_g\in A\star_\af\G$. By \eqref{varphi}, $\psi(a_g\delta_g)\in (C\star_\gm\G(x))_{g_x}$. Then $\psi(a_g\delta_g)\in \psi(a_g\delta_g)(C\star_\gm\G(x))\psi(a_g\delta_g)$ which implies that $a_g\de_g\in (a_g\delta_g)(A\star_\af\G) (a_g\delta_g)$. Therefore  $A\star_\af\G$ is $\G$-graded von Neumann regular.} %{\color{red} If $B$ is $\G$-graded von Neumann regular then $B$ is $\G(x)$-graded von Neumann regular? Para responder isso temos que olhar como sao os ideais da acao parcial $\gamma$. Acho que nao vale a pena, mas se alguem se animar... }
\end{remark}

The next auxiliary result has immediate proof which will be omitted.

\begin{lem1}\label{vnr} Let $\{I_j\}_{1\leq j\leq n}$ be a family of ideals of a ring $R$ that satisfies $I_iI_j=I_jI_i=0$ for all $i\neq j$. Then $\sum_{j=1}^n I_j$ is von Neumman regular if and only if  $I_j$ is von Neumman regular, for all $1\leq j\leq n$.
\end{lem1}

Now we present an extension of Theorem 3 of  \cite{HY} for the context of groupoid graded rings.
\begin{prop1}\label{gvng}  Let $R=\bigoplus_{g\in \G}R_g$ be a  $\G$-graded  ring and $R_0=\bigoplus_{z\in \G_0}R_z$. Then the following assertions hold.
	\begin{enumerate}[\rm (i)]
		\item If $R$ is graded von Neumann regular then  $R_0$ is von Neumann regular.
		\item  If $R_0$ is von Neumann regular, $R$ is strongly $\G$-graded and unital, then $R$ is graded von Neumann regular.
	\end{enumerate}
\end{prop1}

\begin{proof} (i) By Lemma \ref{vnr}, it is enough to show that $R_z$ is von Neumann regular, for all $z\in \G_0$. Let $a\in R_z$. Then there is $r\in R$ such that $a=ara$. Assume that $r=\sum_{g\in \G}r_g$, with $r_g\in R_g$ for all $g\in \G$. Notice that $R_zR_z\subset R_z$. Hence $ara-ar_za=a-ar_za\in R_z$ and consequently
	$$ara-ar_za=\sum_{g\in \G\setminus\{z\}}ar_ga \in R_z\cap \left( \sum\limits_{g\in \G\setminus\{z\}}R_g\right)=0.$$ 
	Thus $a=ara=ar_za\in aR_za$ and whence $R_z$ is von Neumann regular. \smallbreak\smallbreak
	
	\noindent (ii) Let $a\in R_g$ and $g\in \G$. Since $R$ is strongly $G$-graded, ${R_g}$ is a left $R_{t(g)}$-module projective and finitely generated (see the proof of Lemma \ref{noet}). Now,  $R_{g\m}a\subseteq R_{s(g)}$ and $R_{s(g)}R_{g\m}a=R_{g\m}a$,  that is, $R_{g\m}a$ is a left ideal of $R_{s(g)}$. 
	%Consider a set of generators  $\{y_1, y_2, \cdots, y_l\}$  of $R_{g\m}$ as a left $R_{s(g)}$-module. Then 
	%\[R_{g\m}a=\sum_{i=1}^lR_{s(g)}y_ia=\sum_{i=1}^lR_{s(g)}z_i,\quad\text{where }\, z_i=y_ia\in R_{s(g)},\quad 1\leq i\leq l. \]
	%Thus $R_{g\m}a$ is finitely generated as a left $R_{s(g)}$-module. 
	Since $R_{s(g)}$ is von Neumann regular by Lemma \ref{vnr},  there is an idempotent $u\in R_{s(g)}$ such that $R_{g\m}a=R_{s(g)}u$. Note also that, by  Proposition  2.1.1 of \cite{lundstrom2004}, $1_R=\sum_{z\in \G_0}1_z$ and  $R_{s(g)}$ is a unital ring with identity  element $1_{s(g)}$. Consequently    $u=1_{s(g)}u=ra$, for some $r\in R_{g\m}$. Moreover $R_{t(g)}a=R_gR_{g\m}a=R_gR_{s(g)}u=R_gu$ which implies that $a=1_{t(g)}a=bu$ for some $b\in R_g$. Then $ara=au=bu^2=bu=a$, that is, $R$ is graded von Neumann regular.
\end{proof}

{\it For the rest of this subsection we will assume that $\G$ is finite}. Recall that the trace map ${\rm tr}_\alpha:A\to A$ associated to the partial action $\alpha$ is defined by $${\rm tr}_\alpha(a)=\sum_{g\in \G}\alpha_g(a1_{g\m}),\quad a\in A.$$
Since any group is a groupoid, the trace map is defined to partial group actions. Notice that ${\rm tr}_\alpha(1_A)=\sum_{g\in \G}1_{g}$. The trace map has a main role in the Galois theory for partial groupoid actions. More details can be seen in \cite{BP}. \smallbreak

In the sequel we prove two auxiliary results.

\begin{lem1}\label{lem-jacobson1} Let $(C,\gamma)$ be the partial action of $\G(x)$ on $C=A\ast_{\beta}\G_0^2$ given by \eqref{Ch} and \eqref{gamma-h} and $a\in A_x$. Then:
	\begin{enumerate}[\rm (i)]
		\item ${\rm tr}_{\gamma}(1_C)=\sum_{z\in \G_0}\af_{\tau_z}({\rm tr}_{\af_{\G(x)}}(1_x))\delta_{(z,z)}$. \smallbreak
		\item An element $b=\sum_{z\in \G_0}\af_{\ta_z}(a)\delta_{(z,z)}\in C$ is invertible in $C$ if and only if $a$ is invertible in $A_x$. In this case, $b^{-1}=\sum_{z\in \G_0}\af_{\ta_z}(a^{-1})\delta_{(z,z)}$. 
	\end{enumerate}
\end{lem1}
\begin{proof}
	Since $C_h=C1'_{h}$ with $1'_h=\sum_{z\in \G_0}\af_{\ta_z}(1_h)\delta_{(z,z)}$ and $1_C=\sum_{y\in\G_0}1_y\delta_{(y,y)}$ we have  
	\begin{align*}
		{\rm tr}_\gamma(1_C)&=\sum_{h\in \G(x)}\sum_{z\in \G_0}\af_{\tau_z}(1_{h})\delta_{(z,z)}=\sum_{z\in \G_0}\af_{\tau_z}\Big(\sum_{h\in \G(x)}1_{h}\Big)\delta_{(z,z)}\\
		&=\sum_{z\in \G_0}\af_{\tau_z}\Big(\sum_{h\in \G(x)}\af_{\G(x)}(1_x1_{h\m})\Big)\delta_{(z,z)}=\sum_{z\in \G_0}\af_{\tau_z}\big({\rm tr}_{\af_{\G(x)}}(1_x)\big)\delta_{(z,z)},
	\end{align*}	
	and whence (i) follows. For (ii), assume that $a$ is invertible in $A_x$ and consider the element $c=\sum_{z\in \G_0}\af_{\ta_z}(a^{-1})\delta_{(z,z)}\in C$. Thus, 
	$bc=\sum_{z\in \G_0}\af_{\tau_z}(1_x)\delta_{(z,z)}=\sum_{z\in \G_0}1_z\delta_{(z,z)}=1_C$ and whence $c=b^{-1}$. Conversely, consider $c\in C$ the inverse of $b$. Then, there are elements $a_u\in B_u=A_{t(u)}$ such that $c=\sum_{u\in \G_0^2}a_u\delta_u$. As $\af_{\ta_{t(u)}}:A_x\to A_{t(u)}$ is a ring isomorphism, for each $u\in \G_0^2$, there is $a_{x,u}\in A_x$ such that $a_u=\af_{\ta_{t(u)}}(a_{x,u})$. Hence,
	\begin{align*}
		\sum_{z\in \G_0}1_z\delta_{(z,z)}=1_C=bc&=\sum_{z\in \G_0}\sum_{u\in \G(\,\,,\,z)}\af_{\ta_z}(aa_{x,u})\delta_u\\&=\sum_{z\in \G_0}\sum_{y\in \G_0}\af_{\ta_z}(aa_{x,(y,z)})\delta_{(y,z)}.
	\end{align*}
	Thus, $a_{x,(y,z)}=0$ if $y\neq z$ and $aa_{x,(z,z)}=1_x$. Then $a_{x,(z,z)}=a^{-1}$ and we have that $c=\sum_{z\in \G_0}\af_{\ta_z}(a^{-1})\delta_{(z,z)}$. 
\end{proof}

\begin{lem1}\label{lem-jacobson2} Let $(A,\beta)$ be the global action of $\G_0^2$ on $A$ given in the previous section. Then ${\rm tr}_\beta(1_A)=\vert \G_0 \vert 1_A$.
\end{lem1}
\begin{proof}
	It is clear that 
	\begin{align*}
		{\rm tr}_\beta(1_A)&= \sum\limits_{u\in \G_0^2 }\beta_u(1_{s(u)})=\sum\limits_{u\in \G_0^2 }1_{t(u)}= \sum\limits_{(y,z)\in \G_0^2 }1_{z}=\vert \G_0 \vert\sum\limits_{z\in \G_0 }1_{z}=\vert \G_0 \vert 1_A,
	\end{align*}	
	because $A=\bigoplus\limits_{z \in \G_0}A_z$.
\end{proof}

Now we can prove the  main result of this subsection.\smallbreak

\begin{teo1}\label{vn1}
	Assume that $\G$ is finite  and let $\beta$ be the global action of $\G$ on $A$ given in \eqref{action-beta}. Then the following assertions hold.
	\begin{enumerate}[\rm (i)]
		\item  If $A$ is von Neumann regular and $I$ is a left principal ideal of $A\star_{\bt}\G_0^2$, then $I$ is an $A$-direct summand of $A\star_{\bt}\G_0^2$. \smallbreak
		
		\item If  $|\G_0|1_A$ is invertible in $A$ then the following statements are equivalent:\smallbreak
		\begin{enumerate} [\rm (a)]
			\item $A$ is von Neumann regular,\smallbreak
			\item $A\star_{\bt}\G_0^2$ is graded  von Neumann regular,\smallbreak
			\item $A\star_{\bt}\G_0^2$ is von Neumann regular,\smallbreak
			\item  $A\star_\af\G$ is  graded von Neumann regular.\smallbreak
		\end{enumerate}
		Moreover, under the assumption that at least one of the above statements holds and  ${\rm tr}_{\af_{\G(x)}}(1_x)$ is invertible in $A_x$ the following additional statement  holds: \smallbreak
		\begin{enumerate} [\rm (e)]
			\item $A\star_\af\G$ is von Neumann regular.
		\end{enumerate}
	\end{enumerate}
\end{teo1}

\begin{proof} For (i), consider $u,v\in \G_0^2$ such that $s(v)=t(u)$ and $a\in B_u=A_{t(u)}$. Then $A_{t(v)}\delta_v\cdot a\delta_u=A_{t(v)}\beta_v(a1_{s(v)})\delta_{vu}$. Suppose that $I=(A\star_{\bt}\G_0^2)l$, with $l=\sum_{u\in \G_0^2}a_u\delta_u$. Then
	$$I=\sum_{u,v\in \G_0^2}A_{t(v)}\delta_va_u\delta_u=\sum_{t(u)=s(v)}A_{t(v)}\beta_v(a_u1_{s(v)})\delta_{vu}.$$
	For each $u\in \G_0^2$, let $$I_u:=\sum_{s(v)=t(u)} A_{t(v)}\beta_v(a_u1_{s(v)})=\sum_{s(v)=t(u)} A1_{t(v)}\beta_v(a_u1_{s(v)}).$$ It is clear that $I_u$ is a finitely generated left ideal of $A$. Since $A$ is von Neumann, $I_{u}$ is a direct summand of the left $A$-module $A$. Hence $I$ is  an $A$-direct summand of $C:=A\star_{\bt}\G_0^2$.\smallbreak

	For (ii),  using the identification map (\ref{inclusion-map}) we have by  \eqref{isosos} that $C$ is strongly graded and $(A\star_{\bt}\G_0^2)_0=C_0=A$. Thus, the equivalence  (a) $\Leftrightarrow $ (b) follows from Proposition \ref{gvng}. In order to prove (a) $\Rightarrow $ (c), consider a principal left ideal $I$ of $C$. By (i), $I$ is an $A$-direct summand of $C$. Also, by Lemma \ref{lem-jacobson2}, $\operatorname{tr}_{\beta}(1_A)=\vert \G_0 \vert1_A$. Then, Theorem 5.5 of \cite{FMP} implies that $I$ is a $C$-direct summand of $C$. Thus $C$ is von Neumann regular.
	For (c) $\Rightarrow $ (d), observe that $C\star_\gm\G(x)$ is an epsilon strongly  $G(x)$-graded ring such that  $(C\star_\gm\G(x))_x=C$. Then,  from Corollary 3.11 of \cite{La2} we obtain that  $C\star_\gm\G(x)$ is $G(x)$-graded von Neumann regular.  By Remark \ref{obs-g-grade},  $A\star_\af\G$ is $\G$-graded von Neumann regular. Finally, note that (d) $\Rightarrow $ (a) follows from Proposition \ref{gvng}  because $(A\star_\af\G)_0=\bigoplus\limits_{z\in \G_0}A_z\delta_z\simeq A$.  To finish the proof we suppose also that  ${\rm tr}_{\af_{\G(x)}}(1_x)$ is invertible in $A_x$ and prove (e). Indeed, it follows from itens (i) and (ii) of Lemma \ref{lem-jacobson1} that ${\rm tr}_{\gamma}(1_C)$ is an invertible element in  the von Neumann regular ring $A*_{\beta}\G_0^2$.   Hence, Theorem 4.3 of \cite{FL} implies that $(A\star_{\bt}\G_0^2)\star_ \gamma\G(x)$ is von Neumann regular. Thus, we obtain from Theorem \ref{teo-decomp} that $A\star_\af\G$ is von Neumann regular. 
\end{proof}

\begin{remark}{\rm Notice that the assumption  $|\G_0|1_A $ is invertible in $A$  in Theorem \ref{vn1} was only used to prove (a) $\Rightarrow $ (c).}
\end{remark}
\subsection{Jacobson radical} Here we investigate the relation between $J(A)$ and $J(A*_\alpha \G)$. 
The following properties of $J(A)$ will be useful later.

\begin{lem1} \label{jsum}The following statements hold.
	\begin{enumerate}[ \rm (i)]
		\item $J(A)$ is an $\alpha$-invariant ideal of $A$, that is, $\alpha_g(J(A)\cap A_{g\m})\subset J(A)$ for all $g\in \G$.\smallbreak
		\item $J(A_{t(g)})1_g=J(A_g),$ for all $g\in \G$; in particular $J(A_g)$ is unital.\smallbreak
		\item $J(A)=\bigoplus_{z\in \G_0} J( A_z)$.\smallbreak
		\item $\U(A)=\U(A\star_\alpha \G) \cap A$. \smallbreak
		\item $J(A)  = A\cap J(A\star_\alpha \G)$.
	\end{enumerate}
\end{lem1}

\begin{proof} 
	Let $g\in \G$. As $J(A)$ is an hereditary radical and $\af_g$ is a ring isomorphism one has
	$$\alpha_g(J(A)\cap A_{g^{-1}})=\alpha_g(J(A_{g^{-1}}))=J(\alpha_g(A_{g^{-1}}))=J(A_g)=J(A)\cap A_g,$$
	and whence (i) follows. Notice that (iii) is immediate because $A=\bigoplus_{z\in \G_0}A_z$. To prove (ii), observe that $A_{t(g)}=A_g\oplus A_{t(g)}(1_{t(g)}-1_g)$. Consequently, we have that $J(A_{t(g)})=J(A_g)\oplus J\big(A_{t(g)}(1_{t(g)}-1_g)\big)$ which implies that $J(A_{t(g)})1_g=J(A_g)$.\smallbreak
	
	The inclusion $\subseteq$ in (iv) is clear. For the  reverse, consider the $(A,A)$-bimodule epimorphism
	\begin{equation}\label{pia}\pi_A: A\star_\alpha \G\to A, \quad \pi_A\Big(\sum_{g\in \G} a_g\delta_g\Big)=\sum_{e\in \G_0} a_e\delta_e.
	\end{equation} 
	Let $u\in \U(A\star_\alpha \G) \cap A$. Then there is $v\in A\star_\alpha \G$ such that $uv=vu=1_{A\star_\alpha \G}$. Since $\pi_A(u)=u$ it follows that 
	$u\pi_A(v)=\pi_A(v)u=1_{A}$ and consequently $u\in \U(A)$.
	
	%The proof of the item (v) is inspired in the proof of Theorem 2.5 of \cite{P}. 
	(v) Let $a=\sum_{g\in \G}a_g\delta_g\in A\star_\alpha \G$, where $a_g\in A_g$, for all $g\in \G$. Observe that  $a_g\delta_g=(1_g\delta_g)\big(\af_{g\m}(a)\delta_{s(g)}\big)=(1_g\delta_g)\varphi(\af_{g\m}(a))\in (1_g\delta_g)A,$ being $\varphi$  the inclusion map defined in \eqref{inclusion-map}. Hence $a\in \sum_{g\in \G}(1_g\delta_g)A$ and  $A\star_\alpha \G=\sum_{g\in \G}(1_g\delta_g)A$. Let $M$ be a simple right $A\star_\alpha \G$-module. Then, given a non-zero element $m\in M$, we have $M=m(A\star_\alpha \G)=\sum_{g\in \G}\left(m(1_g\delta_g)\right)A$, thus $M$ is a finitely generated $A$-module.  Fix $N:=MJ(A)$. We shall check that $N=\{0\}$. If $g\in \G$ and $t\in A$, then 
	\begin{align*}
		N(t1_g\delta_g)&=M(J(A)t)(1_g\delta_g)\subseteq MJ(A)(1_g\delta_g)\stackrel{{\rm (iii)}}{=}MJ(A_{t(g)})(1_g\delta_g)\\&\stackrel{{\rm (ii)}}{=}MJ(A_{g})(1_g\delta_g)
		\stackrel{(\ast)}{=}M(1_g\delta_g)J(A_{g\m})\subseteq M(1_g\delta_g)J(A_{s(g)})\\&\subseteq MJ(A)=N.
	\end{align*}
	Note that $(\ast)$ is true because $b(1_g\delta_g)=(1_g\delta_g)(\af_{g\m}(b)\delta_{s(g)})$ and $\af_{g\m}(b)\in J(A_{g\m})$, for all $b\in J(A_g)$. Therefore $N(t1_g\delta_g)=N,$ for all $t\in A$ and $g \in \G$, which implies that $N$ is a right $A\star_\alpha \G$-submodule of  $M$. Since $M$ is simple, $N= M$ or $N=\{0\}$.  If $N=M$ then $MJ(A)=M$. Using that $M$ is an $A$-module finitely generated, we obtain from Lemma 2.4 of \cite{P} that $M=\{0\}$ which is a contradiction. Therefore $N=\{0\}$ and thus $J(A)$ annihilates every simple $A\star_\alpha \G$-module. Consequently $J(A)\subseteq J(A\star_\alpha \G)$. For the reverse, we recall that $a\in J(A)$ if and only if $1-xa\in \mathcal{U}(A)$. Hence it follows from (iv) that $A\cap J(A\star_\alpha \G)\subset J(A)$. 
\end{proof}

We end this subsection with the following.

\begin{teo1} \label{jaco}If ${\rm tr}_{\af_{\G(x)}}(1_x)$ is invertible in $A_x$ and $\vert G_0\vert1_A$ is invertible in $A$ then $$J(A*_\alpha \G)\simeq J(A)*_\alpha \G.$$
\end{teo1}
\begin{proof}  It is immediate from Theorem \ref{teo-decomp} that $J(A*_\alpha \G)\simeq J\big((A\star_{\bt}\G_0^2)\star_ \gamma\G(x)\big)$.
	By Lemma \ref{lem-jacobson1}, ${\rm tr}_{\gamma}(1_C)$ is invertible in $C=A\star_{\bt}\G_0^2$. Then, it follows from Proposition 6.7 of \cite{FL} that $J((A\star_{\bt}\G_0^2)\star_ \gamma\G(x))=J(A\star_{\bt}\G_0^2)\star_ \gamma\G(x)$. Also, by Lemma \ref{lem-jacobson2}
	and Theorem 5.7 of \cite{FMP} we have that $J(A*_\beta \G^2_0)=J(A)*_\beta \G^2_0$. Hence  $$J((A\star_{\bt}\G_0^2)\star_ \gamma\G(x))=(J(A)\star_{\bt}\G_0^2)\star_ \gamma\G(x).$$
	By Lemma \ref{jsum} (i) and Remark \ref{ob-invar-ideal} the restriction of $\alpha$ to $J(A)$, which will be denoted by $\af|_{J(A)}$, is also a group-type partial action. Moreover, Lemma \ref{jsum} (ii) implies that   $\af|_{J(A)}$ is unital. Applying Theorem \ref{teo-decomp} for the partial groupoid action $\af|_{J(A)}$ we obtain that $J(A)*_\alpha \G\simeq(J(A)\star_{\bt}\G_0^2)\star_ \gamma\G(x)$. Consequently $J(A*_\alpha \G)\simeq J(A)*_\alpha \G,$ as desired.
\end{proof}

\subsection{Right perfect property} 
A ring $R$ is said to be  {\it right perfect} if every right $R$-module has a projective cover. It is known that $R$ is right perfect if and only if $R/J(R)$ is right artinian and $J(R)$ is right $T$-nilpotent (that is, for every infinite sequence of elements of $J(R)$, there is an $n$ such that the product of first $n$ terms is zero). Also, $R$ is right perfect if and only if $R$ satisfies the descending chain condition (DCC) on principal left ideals. These equivalences for the notion of right perfect ring can be seen, for instance, in \cite{AF}. 
\smallbreak

The next result has the same proof of Lemma 3.4 of \cite{P}. 

\begin{lem1} \label{herej2} Let $R$ be a ring and $\theta$ a unital global action of $\G$ on $R$. Assume that $I$ is a $\theta$-invariant ideal of $R.$ Then the following assertions hold.
	\begin{enumerate}[ \rm (i)]
		\item If $I$ is nilpotent, then $I\star_\theta \G$ is nilpotent. \smallbreak
		\item If $I$ is  right $T$-nilpotent, then $I\star_\theta \G$ is  right $T$-nilpotent.
	\end{enumerate}
\end{lem1}

We need one more auxiliary result. \smallbreak

\begin{lem1}\label{quasi} Let $\theta=(S_g, \theta_g)$ be a  global action of a groupoid $\G$ on a ring $S$ and $I$ a $\theta$-invariant ideal of $S.$ Then $\theta$ induces a global action $\overline{\theta}$ of $\G$ on $S/I$. Moreover, the map
	\[\Psi:\big(S\star_\theta\G\big)/\big(I\star_\theta\G\big)\to\big(S/I\big)\star_{\overline{\theta}}\G,\quad a\delta_g+(I\star_\theta\G)\mapsto (a+I)\delta_g, \quad a\in S_{s(g)}, \]
	is a well-defined ring isomorphism. 
\end{lem1}
\begin{proof} Let $g\in \G$. Put  $B_g=B_{t(g)}:= \big(S_{t(g)}+I\big)/{I}$. Clearly, $B_g=B_{t(g)}$ is an ideal of $S/I$. As in the proof of Lemma 2.2 of \cite{FL}, the ring isomorphism $\theta_g:S_{g\m}\to S_g$ induces the following ring isomorphism
	\begin{align*}
		&\overline{\theta}_g:B_{g\m} \to B_{g},&&\overline{\theta}_g(a+I)=\theta_g(a)+I,&a\in S_{s(g)}.& 
	\end{align*} 
	Since $I$ is $\theta$-invariant it follows that $\overline{\theta}_g$ is well-defined ring isomorphism, for each $g\in \G$.
	It is straightforward to check that $\overline{\theta}=(B_g, \overline{\theta}_g)_{g\in \G}$ is a global action of $\G$ on $S/I$. Finally, consider the surjective ring homomorphism $\psi:S\star_\theta\G\to \big(S/I\big)\star_{\overline{\theta}}\G$ given by $\psi(a\delta_g)=(a+I)\delta_g$, for all $a\in S_{s(g)}$. It is clear that $\ker \psi=I\star_\theta\G$ and whence we obtain the ring isomorphism $\Psi$.  
\end{proof}

The main result of this subsection is the following.

\begin{teo1} \label{rightp}Suppose that $\alpha$ is finite-type. Then $A\star_\af\G$ is right perfect if and only if $A$ is right perfect and $\G$ is finite. 
\end{teo1}
\begin{proof} Suppose that $A\star_\af\G$ is right perfect. Observe that by \eqref{inclusion-map}, $A$ is a subring of $A\star_\af\G$. Moreover $A$ is a right $A$-direct summand of $A\star_\af\G$. From Proposition 2.1 of \cite{PA} follows that $(A\star_\af\G)I\cap A=I$ for every left ideal $I$ of $A$.  Since $A\star_\af\G$ satisfies the (DCC) for left principal ideals follows that $A$ satisfies the (DCC) for left principal ideals. Thus, $A$ is right perfect. 
	Also,  by Theorem \ref{teo-decomp} we have that  $(A\star_{\bt}\G_0^2)\star_ \gamma\G(x)$ is right perfect, where $\gamma$ is finite-type (by Lemma \ref{latoga}). Then, Theorem 34 of \cite{BC} implies that $\G(x)$ is finite. Therefore, by Proposition \ref{group:connec} we conclude that $\G$ is finite. \smallbreak
	
	Conversely, suppose that  $A$ is right perfect and $\G$ is finite. We claim that $C:=A\star_{\bt}\G_0^2$ is right perfect. In fact, fix $B:=A/J(A)$. Notice that $\G_0^2$ acts globally on $A$ by $\beta$ and $J(A)$ is $\beta$-invariant. Hence, by Lemma \ref{quasi}, $\beta$ induces a global action $\overline \beta$ on $B$. Also, $B$ is right artinian because it is semisimple. As $\G_0^2$ is finite we have by Theorem 1.2 of \cite{NOP} that $B*_{\overline\beta}\G_0^2$ is right artinian.  Using Lemma \ref{quasi} we obtain that 
	\begin{equation}\label{firstiso}
		B*_{\overline\beta}\G_0^2\simeq C/\big(J(A)\star_{\bt}\G_0^2\big) %Lemma \ref{jsum} (v) 
	\end{equation} and  $C/\big(J(A)\star_{\bt}\G_0^2\big)$ is right artinian. We claim that $C/J(C)$ is a ring epimorphic image of $C/\big(J(A)\star_{\bt}\G_0^2\big)$. In fact,  for $u\in \G_0^2$ and $a_u\delta_u\in J(A)\star_{\bt}\G_0^2$  we have  $a_u\in J(A)\cap A_{t(u)}$. From Lemma \ref{jsum} (v) follows that $a_u1_{C}\in J(C)$.  Since $J(C)$ is an ideal of $C$  we obtain that $a_u\delta_u=(a_u1_C)(1_u\delta_u)\in J(C)$ and whence $J(A)\star_{\bt}\G_0^2\subset J(C)$. Consequently, the claim is proved and  $C/J(C)$ is right artinian. Now we  show that $J(C)$ is  right $T$-nilpotent. Let $\pi$ be the canonical homomorphism  from $C$ onto $C/\big(J(A)\star_{\bt}\G_0^2\big)$. Notice that
	\begin{equation*} \pi(J(C))=J(C)/\big(J(A)\star_{\bt}\G_0^2\big)\subseteq  J\Big(C/\big(J(A)\star_{\bt}\G_0^2\big)\Big).
	\end{equation*}  
	On the other hand, \eqref{firstiso} implies that  $J(B*_{\overline\beta}\G_0^2)\simeq J\Big(C/\big(J(A)\star_{\bt}\G_0^2\big)\Big),$ consequently $J(C)/\big(J(A)\star_{\bt}\G_0^2\big)$ is nilpotent because it is embedded in the nilpotent ring $J(B*_{\overline\beta}\G_0^2)$. Also, using that $A$ is right $T$-nilpotent, it follows from Lemma \ref{herej2} (ii)  that  $J(A)\star_{\bt}\G_0^2$ is right $T$-nilpotent. Therefore  $J(A\star_{\bt}\G_0^2)$ is right $T$-nilpotent. Hence, $C$ is right perfect. Then by Theorem \ref{teo-decomp} and Proposition 2.1 of \cite{PA} we obtain that $A*_{\alpha} \G$ is right perfect, as desired.
\end{proof}

\subsection{The semiprimary property} 
Recall that a ring $A$ is called {\it semiprimary} if  $A/J(A)$ is  right artinian  and $J(A)$ is nilpotent. 
Clearly every semiprimary ring is right perfect.

\begin{teo1}\label{spp} Suppose that $\alpha$ is finite-type.  Then $A\star_\af\G$ is semiprimary if and only if $A$ is semiprimary and $\G$ is finite.   
\end{teo1}
\begin{proof}  Assume that $A\star_\af\G$ is semiprimary. Using Theorem \ref{teo-decomp} and Lemma \ref{latoga}, it follows from  Theorem 35 of \cite{BC} that $A\star_{\bt}\G_0^2$ is semiprimary and $\G(x)$ is finite. Then $\G$ is finite thanks to Proposition \ref{group:connec}. Since $A\star_{\bt}\G_0^2$ is semiprimary  it is right perfect. Hence, Theorem \ref{rightp} implies that $A$ is right perfect and whence $A/J(A)$ is right artinian. Also, using (v) of Lemma \ref{jsum}  we conclude that $J(A)\subseteq J(A\star_{\bt}\G_0^2)$. Consequently $J(A)$ is nilpotent. \smallbreak
	
	Conversely, suppose that $A$ is semiprimary and $\G$ is finite. Then $A$  is right perfect. As we saw in the proof of Theorem \ref{rightp}, $C:=A\star_{\bt}\G_0^2$ is right perfect. Thus $C/J(C)$ is right artinian.  Since $J(A)$ is nilpotent it follows from  Lemma \ref{herej2} (i) that  $J(A)\star_{\bt}\G_0^2$ is nilpotent. Moreover, as in the proof of Theorem \ref{rightp}, $J(C)/\big(J(A)\star_{\bt}\G_0^2\big)$ is nilpotent. Then $J(C)$ is nilpotent  and whence $C$ is semiprimary. Finally, since $\G(x)$ is finite, the same argument shows that $C\star_ \gamma\G(x)$ is semiprimary. Thus, the result follows from Theorem \ref{teo-decomp}.
\end{proof}

\subsection{Krull dimension}
Let $R$ be a unital ring and $M$ a right $R$-module. The Krull dimension of $M$,  usually denoted by $\operatorname{Kdim} M_R$
is the deviation of the lattice $\mathcal{L}_R(M)$ of $R$-submodules of $M$. The right Krull dimension of a ring $R$ is defined as $\operatorname{Kdim} R:=\operatorname{Kdim} R_R$. 
It is well-known  that the right modules of Krull dimension $0$ are the nonzero right artinian modules and that every right noetherian module has Krull dimension.

\begin{teo1}\label{krull}
	Assume that $\G$ is finite. If $A$ has right Krull dimension then $A\ast_{\af}\G$ has right Krull dimension and $\operatorname{Kdim} A\ast_{\af}\G\leq \operatorname{Kdim} A$.
\end{teo1}

\begin{proof}
	Let $u\in \G_0^2$. Since $A$ has Krull dimension, it follows that $B_u=A_{t(u)}$ has Krull dimension as a right $A$-module. Using that $A_u$ and $A_{t(u)}\delta_{u}$ are isomorphic as right $A$-modules we conclude that $A_{t(u)}\delta_{u}$ has Krull dimension. Since $A\star_{\bt}\G_0^2=\bigoplus_{u\in \G_0^2}A_{t(u)}\delta_{u}$,  Lemma 6.1.14 of \cite{MR} implies that  $C:=A\star_{\bt}\G_0^2$ has Krull dimension as a right $A$-module. Moreover, it is easy to see that $\operatorname{Kdim} C_A=\operatorname{Kdim}A_A=\operatorname{Kdim}A$.
	As in the proof of Proposition 41 of \cite{BC}, we have that
	$$\mathcal{L}_{C}(C)\to \mathcal{L}_{A}(C),\quad N\mapsto N,$$
	is strictly increasing. Thus, $C$ has Krull dimension and $\operatorname{Kdim} C\leq \operatorname{Kdim} C_A= \operatorname{Kdim}A$.
	From Theorem \ref{teo-decomp} and Proposition 41 of \cite{BC} we obtain that $A\star_\af\G$ has right Krull dimension and $\operatorname{Kdim}A\star_\af\G\leq \operatorname{Kdim}C\leq  \operatorname{Kdim}A$.
\end{proof}

\section{Applications to Leavitt path algebras}\label{Leavitt}
The Leavitt path algebra associated to a directed graph  was introduced in \cite{AA} and \cite{AAP} as follows. Denote a directed graph by $E=(E^0,E^1,r,d)$, where $E^0$ and $E^1$ are countable sets and $r,d : E^1 \to E^0$ are maps. The elements of $E^0$ are called \emph{vertices} and the elements of $E^1$ are called \emph{(real) edges}.  
A \emph{sink} is a vertex such that no edge emerges out of it. A \emph{path $\mu$ of length n} in 
$E$ is a sequence of edges $\mu=\mu_1 \cdots \mu_n$ such that $r(\mu_i)=d(\mu_{i+1})$ for $i\in \{1,\ldots,n-1\}$. In such a case, we will write $|\mu|=n,$  $d(\mu):=d(\mu_1)$ is the \emph{source} of $\mu$ and $r(\mu):=r(\mu_n)$ is the \emph{target} of $\mu,$ respectively, a vertex is considered as a path of length $0$. We are using $r,d$, instead $t,s$ to denote the target and source maps of a path because we want to avoid confusion with the notation used for the target and source maps in groupoids. We shall denote $E^{1}(v,\,):=\{ f\in E^1\,:\, d(f)=v \}$.

\begin{def1}{\rm
		Let $E=(E^0,E^1,r,d)$ be a directed graph and $\Bbbk$ a field. The \emph{Leavitt path $\Bbbk$-algebra $L_\Bbbk(E)$ of $E$ with coefficients in $\Bbbk$} is the free associative
		$\Bbbk$-algebra generated by the set $\{v,f,f^{\ast}\, :\, v\in E^0, f\in E^1 \}$ with the following relations:
		\begin{enumerate}[\rm (i)]
			\item for all $v,w\in E^0$, $v^2=v$ and $vw=0$ if $v\neq w$,\vspace{.2em} 
			\item $d(f)f=fr(f)=f$, for all $f\in E^1$,\vspace{.3em}
			\item $r(f)f^*=f^*d(f)=f^*$, for all $f\in E^1$,\vspace{.3em}
			\item for all $f,f'\in E^1$, $f^*f=r(f)$ and $f^*f'=0$ if $f\neq f'$,\vspace{.2em}
			\item $v=\sum\limits_{ f\in E^{1}(v,\,) } ff^*$, for every $v\in E^0$ such that $d^{-1}(v)$ is non-empty and finite.
	\end{enumerate}}
\end{def1}
The  symbols $f^*
$  for $f \in E^1$ are called {\it ghost edges}, also
condition (v) above is called the {\it Cuntz-Krieger relation}.  Notice that $L_\Bbbk(E)=\bigoplus_{v\in E^0}L_\Bbbk(E)v$
and thus $L_\Bbbk(E)$ is  unital if and only if $E^0$  is finite. In this case, $1_{L_\Bbbk(E)}=\sum_{v\in E^0}v$. \smallbreak

In what follows in this section, $E=(E^0,E^1,r,d)$ denotes a directed graph and $L_\Bbbk(E)$ is the corresponding Leavitt path algebra, where $\Bbbk$ is a field.  Our purpose in this section is to apply the results of the previous section to $L_\Bbbk(E)$.

\subsection{$L_\Bbbk(E)$ as a partial skew groupoid ring}
Leavitt path algebras were realized as partial skew groupoid rings in Theorem 3.11 of  \cite{GY}. For readers convenience, here we recall this realization.

Firstly, the free path groupoid $\G(E)$ associated to $E$ is constructed in the following way. The set of objects  $\G(E)_0$ of $\G(E)$ is $E^0$. We extend the maps $r,d$ from $E^{1}$ to $E^{1}\cup \big(E^{1}\big)^\ast$ by putting: $d(f^*)=r(f)$ and $r(f^*)=d(f)$, for all $f\in E^{1}$. Then, for a finite path $\mu=\mu_1\ldots\mu_n$ in $E$ we put $\mu^{\ast}:=\mu_n^{\ast}\ldots\mu_1^\ast$. Denote by 
$$P:	=\{\xi_1\xi_2\cdots \xi_n\,:\, n\in \N, \,\, \xi_i\in E^1\cup (E^1)^*\cup E^0\,\, \,\text{and}\,\,\, r(\xi_i)=d(\xi_{i+1})\}.$$
An element $\xi_1\ldots\xi_kff^*\xi_{k+1}\ldots\xi_n$ in $P$ can be reduced to $\xi=\xi_1\ldots\xi_k\ldots\xi_n$. Similarly, $\xi$ is the reduction of  $$\xi_1\ldots \xi_kd(\xi_{k+1})\xi_{r+1}\ldots\xi_n=\xi_1\ldots\xi_rr(\xi_k)\xi_{r+1}\ldots\xi_n.$$ The elements of $P$ that cannot be reduced are called irreducible. The set of irreducible elements of $P$ will be denoted by $\operatorname{Irr}(P)$ and the reduction of an element $\mu\in P$ will be denoted by $\operatorname{irr}(\mu)$.  The set of morphisms of the $\G(E)$ is $\operatorname{Irr}(P)$. Given $\mu,\xi\in \operatorname{Irr}(P)$ such that  $r(\mu)=d(\xi)$, the composition is $\mu\cdot\xi=\operatorname{irr}(\mu\xi)$. It is clear that $\xi^{\m}=\xi^\ast$, for all $\xi\in \operatorname{Irr}(P)$.

\begin{remark}\label{obs-free} {\rm Notice that in the directed graph we compose from left to right while in the free path groupoid $\G(E)$ we compose from right to left. Thus, {\it we have that $t(\xi)=d(\xi)$ and $s(\xi)=r(\xi)$, for all $\xi\in \G(E)$}. Precisely, if $\xymatrix{\bullet^{v}\ar@{->}[r]^{\xi}& \bullet^{w}}$ is an edge in $E$ then $d(\xi)=v$ and $r(\xi)=w$. However, $\xi$ as an element of $\G(E)$ has $s(\xi)=w$ and $t(\xi)=v$. Clearly, given $\xi,\eta\in E^{1}$ we have that  $\xi\eta$ exists in $E$ if and only if $\eta\xi$ exists in $\G(E)$.}
	% (ii) It is clear that $\G(E)$ is a connected groupoid  if and only if $E$ is a connected directed graph.} 
\end{remark}

Let $W$ be the set of all finite paths in $E.$  Denote by $W^{\infty}$ the set of all infinite paths in $E$, that is, $\xi\in W^\infty$ if and only if $\xi=\xi_1\xi_2\xi_3\cdots$ with $r(\xi_i)=d(\xi_{i+1})$ and $i\in\N$. Also, consider
\begin{align*}
	&W^{(1)}=\{a\in W\setminus E^0\,:\,\operatorname{irr}(a)=a \}\quad\text{and} \\[.1em]
	&W^{(2)}=\{(a,b)\in \big(W\setminus E^0\big)^2\,:\,r(a)=r(b)\text{ and }\operatorname{irr}(ab\m)=ab\m\}.
\end{align*}
We  define a partial action $\theta$ of $\G(E)$ on the set
\begin{align}\label{set-X}
	X=\{\xi\in W \,:\, r(\xi) \text{ is a sink}\}\cup  W^{\infty},
\end{align}
in the following way.
For each $g\in \G(E)$, let $X_g$ be defined by:
\begin{enumerate}[$\quad \circ$]
	\item  $X_v:=\{\xi\in X\,:\, d(\xi)=v\}$, for all $v\in E^0,$  \smallbreak
	
	\item $X_a:=\{\xi\in X \,:\, \xi_1\xi_2...\xi_{|a|}=a\}$, for all $a\in W^{(1)}$,\smallbreak
	
	\item $X_{a^{-1}}:=\{\xi\in X \,:\,d(\xi)=r(a)\}$, for all $a\in W^{(1)}$,   \smallbreak
	
	\item $X_{ab^{-1}}:=X_a$, for all $(a,b)\in W^{(2)}$,\smallbreak
	
	\item $X_g:=\emptyset$, for all other $g \in \G(E)$.
\end{enumerate}
The bijections $\theta_g$, for  $g\in \G(E)$, are defined by: $\theta_v:=\id_{X_v}$, for all $v\in E^0$,  and 
\begin{align*}
	&\theta_a:X_{a^{-1}}\rightarrow X_a,&  \quad \theta_a(\xi)&=a\xi,\quad a\in \G(E)\setminus E^0& \\[.3em]
	&\theta_{ab^{-1}}:X_{ba^{-1}}\rightarrow X_{ab^{-1}},& \theta_{ab^{-1}}(\xi)&=a\xi_{(|b|+1)}\xi_{(|b|+2)}\ldots,&
\end{align*}
for all $(a,b)\in W^{(2)}$. It is clear that $\theta_a$ and $\theta_{ab^{-1}}$ have inverses maps given respectively by:
\begin{align*}
	\theta_{a^{-1}}(\xi)&=\begin{cases*}
		r(a),&  \text{if }$ r(a)$ \text{ is  a sink},\\
		\xi_{|a|+1}\xi_{|a|+2}\ldots,& \text{ otherwise},  
	\end{cases*} \\[.3em]
	\theta_{ba^{-1}}(\xi)&=b\xi_{(|a|+1)}\xi_{(|a|+2)}\ldots 
\end{align*}

By Proposition 3.7 of \cite{GY}, the family of pairs $\theta=(X_g, \theta_g)_{g\in\G(E)}$ is a partial action of $\G(E)$ on  $X$. 
Following \cite{GY}, denote the algebra of functions from $X$ to $\Bbbk$ by $\mathcal{F}(X)$ and  put $\mathcal{F}(X_g):=\{f\in \mathcal{F}(X)\,:\,f\text{ vanishes outside of }X_g\}$ for each $X_{g}\neq \emptyset$ and $\mathcal{F}(X_g):=\{\text{null function}\}$ for each $X_{g}=\emptyset$. Also, define $\tilde{\theta}_g(f)=f\circ \theta_{g^{-1}}$, for all $ f\in \mathcal{F}(X_{g^{-1}})$. It follows from Proposition 3.8 of \cite{GY} that $\tilde{\theta}=\big(\mathcal{F}(X_g),\tilde{\theta}_g\big)_{g\in \G(E)}$ is a partial action of $\G(E)$ on $\mathcal{F}(X)$.

For each $g\in \G(E)$, $1_g$ denotes the characteristic of ${X_g}$, that is, $1_g(x)=1$ if $x\in X_g$ and $1_g(x)=0$ if $x\notin X_g$. Notice that $1_g$ is the identity element of $\mathcal{F}(X_g)$. Let $D(X)$ be the $\Bbbk$-subalgebra of $\mathcal{F}(X)$ generated by $\{1_g \,:\, g\in \G(E)\}$.  Fix the ideal  $D(X)_g:=1_g D(X)$ of $\mathcal{F}(X_g)$, for each $g\in \G$. Observe that 
$D(X)_g$ is generated (as a subalgebra of $D(X)$) by $\{1_g1_h\,:\, h\in \G(E)\}$.
It is clear that the restriction of $\tilde{\theta}_g$ to $D(X)_{g^{-1}}$ is a bijection onto $D(X)_g$.
We will denote by $\lambda_g$ the restriction of $\tilde{\theta}_g$ to $D(X)_{g^{-1}}$, that is, 
\begin{align}\label{ac-lam}
	\lambda_g:D(X)_{g^{-1}}\to D(X)_{g}, \qquad  \lambda_g=\tilde{\theta}_g|_{D(X)_{g^{-1}}}.
\end{align}
Then $\lambda=(D(X)_g,\lambda_g)_{g\in \G(E)}$ is a partial action of $\G(E)$ on $D(X)$. \smallbreak

Now we can enunciate Theorem 3.11 of \cite{GY}.

\begin{teo1}\label{iso-GY}
	$L_{\Bbbk}(E)$ is isomorphic (as a $\Bbbk$-algebra) to $D(X)\star_{\lambda}\G(E)$.
\end{teo1}

\begin{remark}\label{rem-cone-case} {\rm  
		Let $E$ be a directed graph. From Proposition 1.2.14 of \cite{AAM} follows that $L_{\Bbbk}(E)\simeq \oplus_{i\in I}L_{\Bbbk}(E_i)$, where $E=\sqcup_{i\in I}E_i$ is the decomposition of $E$ in its connected components. For each $E_i$, we denote the set given in \eqref{set-X} by  $X_i$, the free path groupoid by $\G(E_i)$ and the action of $\G(E_i)$ on $D(X_i)$ by $\lambda_i$. Thus, Theorem \ref{iso-GY} implies that $L_{\Bbbk}(E)\simeq \oplus_{i\in I} D(X_i)\star_{\lambda_i}\G(E_i)$.}
\end{remark}

\begin{remark}\label{rem-dist-groups}{\rm
		The Leavitt path algebra $L_{\Bbbk}(E)$ can be realized as a partial skew group ring cf. Proposition 3.2 of  \cite{GR}. In this case, the group that acts partially is the free group $\mathbb{F}$ generated by $E^{1}$. On the other hand, suppose that the partial action $\lambda$ given in \eqref{ac-lam} is group-type, that is, there exist $v\in E^0$ and a transversal $\ta(v)$ for $v$ such that $\lambda$ satisfies \eqref{cond1}. Assume that $E^0$ is finite. From Theorem \ref{teo-decomp} it follows that $L_{\Bbbk}(E)$ can be realized again as a partial skew group ring. In this case, the group acting partially is the isotropy group $\G(E)(v)$ which in general is easier to deal than $\mathbb{F}$. Hence, sometimes it is more convenient to realize $L_{\Bbbk}(E)$ as a partial skew group ring using the partial action $\lambda$; see Remark \ref{rem-advantage}.}
\end{remark}

\subsection{Determining when $\lambda$ is group-type} Let $E$, $D(X)$, $\lambda$ and $\G(E)$ be as in the previous Subsection. In order to apply the results of Section 3 for a Leavitt path algebra $L_{\Bbbk}(E)$, we need to find under what conditions $\lambda$ is a group-type partial action of $\G(E)$ on $D(X)$.  
Thanks to Remark \ref{rem-cone-case}, we will assume from now on in this Subsection that $E$ is connected and $E^0$ is finite.
\begin{lem1} \label{decompo} The following assertions hold.
	\begin{enumerate}[\rm (i)]
		\item
		$D(X)=\bigoplus_{v\in E^0}D(X)_v$. \smallbreak
		\item  $X_v\neq \emptyset,$ for all $v\in E^0$.\smallbreak
		
		\item  For every $v\in E^0$ such that $E^1(v,\,)$ is non-empty and finite set, we have that 
		\[1_v=\sum\limits_{ f\in E^1(v,\,)}  1_f\quad \text{and}\quad D(X)_v=\bigoplus \limits_{ f\in E^1(v,\,) }D(X)_f.\]
	\end{enumerate}
\end{lem1}
\begin{proof} (i) We claim that $D(X)=\sum_{v\in E^0}D(X)_v$. In fact, consider an element $a=1_{g_1} 1_{g_2}\cdots 1_{g_n}\in D(X)$ with $g_1, g_2, \dots, g_n\in \G(E)$. Since $1_g1_{t(g)}=1_g$ for all $g\in \G(E)$, we have that $a=1_{t(g_1)}a$. Thus $a\in D(X)_{t(g_1)}\subseteq\sum_{v\in E^0}D(X)_v$. Using that any element of $D(X)$ is a linear combination of elements $a$ as above, we conclude that $D(X)=\sum_{v\in E^0}D(X)_v$. Also, the sum is direct because $D(X)_v=D(X)1_v$ and $1_{v_i}1_{v_j}=\delta_{i,j}1_{v_i},$ for all $v,v_i,v_j\in E^0$. \smallbreak
	
	\noindent (ii) Let $v\in E^0$. If $v$ is a sink then $v\in X_v$ and $X_v$  is non-empty. Suppose that  $v$ is not a sink. In this case, there exists $\xi_1\in E^1 $ such that $d(\xi_1)=v$. If $r(\xi_1)=v_1$ is a sink then $\xi_1\in X_v$. Otherwise, there exists $\xi_2\in E^1 $ such that $d(\xi_2)=v_1$.  If $r(\xi_2)=v_2$ is a sink then $\xi_1\xi_2\in X_v$. Otherwise, we repeat the process. If this process is finite then we obtain an element $\xi_1\xi_2\cdots \xi_n \in X_v$. If the process is infinite then we have an infinite path $\xi=\xi_1\xi_2\cdots $ for which $d(\xi)=v$, that is, $\xi\in X_v$. \smallbreak
	
	\noindent (iii)  Let $v\in E^0$ such that $E^1(v,\,)$ is non-empty and finite. Given $f,f'\in E^1$, with $f\neq f'$ and $f,f'\in d^{-1}(v)$, we have that $X_{f}\cap X_{f'}=\emptyset$. Hence, in order to show the equality on the left side in (iii) it is enough to verify that $X_v=\bigcup\limits_{ f\in E^1(v,\,) } X_f$. The inclusion $\supseteq$ is immediate. For the reverse, consider $\xi\in X_v$. Then there are $f\in E^1$ and $\eta\in W\cup W^{\infty}$ with $r(f)=d(\eta)$ such that $\xi=f\eta$. Then $\xi\in X_f$ and $d(f)=d(\xi)=v$ which implies the inclusion $\subseteq$. The equality $D(X)_v=\bigoplus \limits_{ f\in E^1(v,\,)}D(X)_f$ follows from the fact that $1_f1_p=\delta_{f,p}1_f$, for all $f,p\in E^1.$ 
\end{proof}

For the next result, we recall that if $\xi$ is an edge in $E$ such that $d(\xi)=v$ and $r(\xi)=w$ then $\xi$ as an element of $\G(E)$ has $s(\xi)=w$ and $t(\xi)=v$; see Remark \ref{obs-free}. We also denote $W^{(-1)}:=\{a\m\,:\,a\in W^{(1)}\}$.

\begin{teo1} \label{gtype} Suppose that $E$ is connected. Then $\lambda$ is a group-type partial action of $\G(E)$ on $D(X)$ if and only if 
	there are a vertex $v\in E^0$ of $E$ and a transversal $\tau(v)=\{\tau_{w}:v\to w\,:\, w\in E^0\}$ for $v$ in $G(E)$ such that $\tau_w\in W^{(-1)}\cup W^{(1)}\cup W^{(2)}$ for all $w\in E^0, w\neq v$ and
	\begin{enumerate}[\quad $\circ$]\vspace*{-.05cm}
		\item if $\tau_w\in W^{(-1)}$ then\vspace*{-.1cm}
		\begin{align}
			\{\xi\in X: \xi_1\xi_2...\xi_{|\tau\m_w|}=\tau\m_w\}=\{\xi\in X : d(\xi)=v\}, \label{eqq1-tg} 
		\end{align}	
		\item if $\tau_w\in W^{(1)}$ then\vspace*{-.1cm}
		\begin{align}
			\{\xi\in X: \xi_1\xi_2...\xi_{|\tau_w|}=\tau_w\}=\{\xi\in X : d(\xi)=w\}, \label{eq1-tg} 
		\end{align}	
		\item if $\tau_w=a_wb^{-1}_w$ with $(a_w,b_w)\in W^{(2)}$ then\vspace*{-.1cm}
		\begin{align}
			\{\xi\in X: \xi_1\xi_2...\xi_{|b_w|}=b_w\}=\{\xi\in X : d(\xi)=v\}	, \label{eq2-tg}
		\end{align}	
		\item if $\tau_w=a_wb^{-1}_w$ with $(a_w,b_w)\in W^{(2)}$ then\vspace*{-.1cm}
		\begin{align}
			\{\xi\in X: \xi_1\xi_2...\xi_{|a_w|}=a_w\}&=\{\xi\in X : d(\xi)=w\}, \label{eq3-tg}
		\end{align}	
		%	\item if $\tau_w\notin W^{(-1)}\cup W^{(1)}\cup W^{(2)}$  then\vspace*{-.1cm}
		%\begin{align}
		%	\{\xi\in X: d(\xi)=w\}&=\emptyset. \label{eqq3-tg}
		%\end{align}	
		
	\end{enumerate}
\end{teo1}
\begin{proof} Notice that $\lambda$ is a group-type partial action of $\G(E)$ on $D(X)$ if and only if there are $v\in E^0$ and a transversal $\tau(v)=\{\tau_{w}: v\to w \,:\, w\in \G(E)_0=E^0\}$ for $v$ such that
	$D(X)_{\tau\m_w}=D(X)_v$  and $D(X)_{\tau_w}=D(X)_w$,  for all $w\in E^0$. Consequently, $\lambda$ is a group-type if and only if $1_{\tau\m_w}=1_v$  and $1_{\tau_w}=1_w$, for all $w\in E^0$, which is  equivalent to 
	\begin{align}\label{xeq}X_{\tau\m_w}=X_v=\{\xi\in X \,:\, d(\xi)=v\},\quad X_{\tau_w}=X_w=\{\xi\in X \,:\, d(\xi)=w\},
	\end{align} 
	for all $w\in E^0$. Using the definition of the ideals $X_g$, $g\in \G(E)$, given in the previous subsection and  (ii) of Lemma \ref{decompo}, we obtain directly from \eqref{xeq} that $\lambda$ is group-type if and only if the transversal $\ta(v)$ satisfies  \eqref{eqq1-tg}, \eqref{eq1-tg}, \eqref{eq2-tg} and \eqref{eq3-tg}.% and \eqref{eqq3-tg}.
\end{proof}

\begin{remark}\label{ssink} {\rm Notice that the inclusions $\subset $ in equations  \eqref{eqq1-tg}, \eqref{eq1-tg}, \eqref{eq2-tg} and \eqref{eq3-tg} always hold. Also, suppose that $v\in E^0$ is a sink and that there exists a transversal $\tau(v)=\{\tau_{w}:v\to w\,:\, w\in E^0\}$ in $\G(E)$ for $v$ contained in $W$. In this case $\ta(v)\setminus\{v\}\subset   W^{(1)}$  and $\lambda$ is a group-type partial action of $\G(E)$ on $D(X)$ with respect to the transversal $\tau(v)$ if and only if  \eqref{eq1-tg} holds.}
\end{remark}

For the reader's convenience we recall that a path $\xi=\xi_1\xi_2\ldots \xi_n$ in $E$ is {\it closed} if $d(\xi)=r(\xi)$. If $\xi$ is closed  and $d(\xi_i)\neq d(\xi_j)$, for every $i\neq j$, then $\xi$ is called {\it a  cycle}. For a cycle $\xi$ in $E$, $\xi^{\infty}$ denotes the infinity path $\xi\xi\ldots$ in $E$. A directed graph is said {\it acyclic} if it has no cycle.\smallbreak

The next result gives necessary and sufficient conditions on $E$ for $\lambda$ to be group-type in the case that $E$ has a sink and $|E^{0}|\geq 2$.
\begin{prop1}\label{daniel} Suppose that $|E^{0}|\geq 2$ and  that $E$ has a sink $v.$ Then $\lambda$ is a group-type partial action of $\G(E)$ on $D(X)$  if and only if $|d^{-1}(w)|= 1$, for all $w\in E^{0}\setminus\{v\}$ and $E$ is acyclic. 
\end{prop1}
\begin{proof}  Let $v\in E^{0}$ a sink. Since $\lambda$ is group-type, by (i) of Remark \ref{obs-not-depend} and Theorem \eqref{gtype}, we can take a transversal $\tau(v)=\{\ta_{w}\}_{w\in E^0}$ in $\G(E)$ for $v$ such that  $\tau_w\in W^{(-1)}\cup W^{(1)}\cup W^{(2)}$ for all $w\in E^0, w\neq v$, and satisfies \eqref{eqq1-tg}, \eqref{eq1-tg}, \eqref{eq2-tg} and \eqref{eq3-tg}. Let $w\in E^0\setminus\{v\}$. If $\tau_w\in W^{(-1)}$ then $\ta\m_w\in W$ and $t(\ta\m_w)=v$. Hence, there exists $\xi\in E^{1}$ such that $d(\xi)=v$ which is a contradiction because $v$ is a sink. If $\ta_w=a_wb\m_w$ with $(a_w,b_w)\in W^{(2)}$ then $t(b_w)=v$. Again, there is $\xi\in E^{1}$ such that $d(\xi)=v$ which is an absurd. Hence, $\ta_w\in W^{(1)}$. In particular $d(\tau_w)=w$ and  $|d^{-1}(w)|\geq 1$.	
	Suppose that there exists $w_0\in E^{0}$ with $|d^{-1}(w_0)|>1$. If $\ta_{w_0}=\xi_1\ldots\xi_r\in W^{(1)}$ then there is $\eta_1\in E^1$ with $\eta_1\neq \xi_1$ and $d(\eta_1)=w_0$. If $r(\eta_1)$ is a sink then $\eta_1\in X$. Otherwise, there is $\eta_2\in E^1$ with $d(\eta_2)=r(\eta_1)$. If $r(\eta_2)$ is a sink then $\eta_1\eta_2\in X$. Otherwise, we repeat the process. If the process is finite then we obtain an element $\eta=\eta_1\ldots\eta_t\in X$ such that $d(\eta)=w_0$. If the process is infinite, then there is an infinite path $\eta=\eta_1\eta_2\ldots\in X$ such that $d(\eta)=w_0$. In both cases, we have that $\eta\in \{\xi\in X : d(\xi)=w_0\} \setminus \{\xi\in X: \xi_1\xi_2...\xi_{|\tau_{w_0}|}=\tau_{w_0}\}$ which implies that \eqref{eq1-tg}  does not hold and we have a contradiction. 
	Now suppose that $\nu=\nu_1\nu_2\ldots \nu_n$ is a cycle in $E$. Thus   
	$\nu^\infty\in X$. If $\nu^\infty=\tau_{d(\nu_1)}\eta$ for some path $\eta$ in $E$ then $v=r(\tau_{d(\nu_1)})=d(\eta),$ which is a contradiction. Hence  $\nu^\infty\in \{\xi\in X : d(\xi)=d(\nu_1)\} \setminus \{\xi\in X: \xi_1\xi_2...\xi_{|\tau_{d(\nu_1)}|}=\tau_{d(\nu_1)}\}$ and again we have a contradiction. So  $E$ is acyclic. 
	
	Conversely, since $E^0$ is finite, $E$ is acyclic and  $|d^{-1}(w)|= 1$, for every $w\in E^{0}\setminus v$, it follows that there are no infinite paths in $E$.  Also, for any $w\in E^0\setminus\{v\}$, there is a unique (finite) path $\tau_w$ in $E$ with $d(\tau_w)=w$ and $r(\tau_w)=v$. Indeed, write $w=w_1$ and let $\eta_1\in E^1$ with $d(\eta_1)=w_1$. As $E$ is acyclic we have $w_2=r(\eta_1)\neq w_1$. If $w_2=v$  we take $\tau_w=\eta_1$ or else  we repeat the process.  Since $E^0$ is finite, the process is also finite. Thus, there are vertices $w_1=w,\ldots,w_n=v$ with $w_i\neq w_j$ for all $1\leq i\neq j\leq n$ and edges $\eta_1,\ldots,\eta_n$ with  $d(\eta_i)=w_i$ and  $r(\eta_i)=w_{i+1}$, for all $1\leq i\leq n-1$. 
	Hence $\tau_w=\eta_1\cdots \eta_n$  satisfies $d(\ta_w)=w$ and $r(\ta_w)=v$. Suppose that $\xi$ is a path in $E$ such that $d(\xi)=w$ and $r(\xi)=v$. Using that $E$ has no infinite paths, we conclude that $\xi=\xi_1\ldots\xi_r$. Since $|d^{-1}(w)|= 1$, we have that $\eta_1=\xi_1$. Hence, $d(\eta_2)=r(\eta_1)=r(\xi_1)=d(\xi_2)$ which implies $\xi_2=\eta_2$.  In this way one obtains $r=n$ and $\eta_i=\xi_i$,  for all $1\leq i\leq n$. Therefore $\xi=\tau_w$. Finally, notice that 	$\{\mu\in X: \mu_1\mu_2...\mu_{|\tau_w|}=\mu_w\}=\{\mu\in X : d(\mu)=w\}=\{\ta_w\}$ and whence  \eqref{eq1-tg} holds. Thus, $\lambda$ is a group-type partial action of $\G(E)$ on $D(X)$ with transversal $\tau(v)=\{\ta_w\}_{w\in E^0}$ for $v$.  
\end{proof}

\begin{cor1}\label{cor-matrices}
	Suppose that $|E^{0}|\geq 2$ and  that $E$ has a sink $v$. If $\lambda$ is a group-type partial action of $\G(E)$ on $D(X)$ then 
	\[L_{\Bbbk}(E)\simeq D(X)\star_{\beta} (E^0)^2\simeq {\rm M}_n(\Bbbk).\]
\end{cor1}
\begin{proof}
	Let $v$ be the unique sink of $E$. By (i) of Remark \ref{obs-not-depend} and Theorem \eqref{gtype}, there exists a transversal $\tau(v)=\{\ta_{w_j}\,:\,1\leq j\leq n\}$ in $\G(E)$ for $v$ that satisfies \eqref{eqq1-tg}, \eqref{eq1-tg}, \eqref{eq2-tg} and \eqref{eq3-tg}.
	Since $\G(E)(v)$ is the trivial group $\{v\}$, it follows from Theorem \ref{teo-decomp} and Theorem \ref{iso-GY} that $L_{\Bbbk}(E)$ is isomorphic  to $D(X)\star_{\beta} (E^0)^2$, where $\beta=(B_u,\beta_u)_{u\in (E^0)^2}$ is the global action of $\G(E)_0^2=(E^0)^2$ on $D(X)$ given by \eqref{action-beta}.
	Explicitly, if $u=(w_i,w_j)\in (E^0)^2$ then $B_u=D(X)_{w_j}=D(X)1_{w_j}$. As we saw in the proof of Proposition \ref{daniel}, $X_{\ta_{w}}=\{\ta_w\}$, for all $w\in E^{0}$. Thus $1_{w_j}1_g=0$ if $\tau_{w_j}\notin X_g$ or $1_{w_j}1_g=1_{w_j}$ if $\tau_{w_j}\in X_g$, for all $g\in \G(E)$. Hence $D(X)1_{w_j}=\Bbbk\langle 1_{w_j}\rangle$ which is isomorphic (as algebra) to $\Bbbk$. Also, notice that $\beta_u:\Bbbk\langle 1_{w_i}\rangle \to \Bbbk\langle 1_{w_j}\rangle$ is given by $\beta_u(a1_{w_i})=a1_{w_j}$, for all $a\in \Bbbk$. Hence it is clear that 
	\begin{align*}
		\psi: D(X)\star_{\beta} (E^0)^2\to{\rm M}_n(\Bbbk),\quad \psi(a1_{w_j}\delta_{(w_i,w_j)})=ae_{ji},\quad a\in\Bbbk,
	\end{align*}
	is an algebra isomorphism, where $e_{ji}$ is the elementary matrix that has $1$ in the $(j,i)$-entry and $0$ in the other entries. 
\end{proof}

\begin{remark}{\rm
		The result of the previous corollary is known; see, for instance, Theorem 2.6.17 of \cite{AAM}.   }
\end{remark}

\begin{exe}\label{toe} {\rm Consider  the \textit{Toeplitz graph} $E_T$ given by 
		
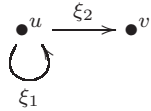
\begin{align*}
			\xymatrix{ &\bullet^{u}\ar@(dl,dr)_{\xi_1} \ar@{->}[r]^{\xi_2}&\bullet^{v}} 
		\end{align*}
		Then follows by Proposition \ref{daniel} that $\lambda$ is not a group-type partial action of $\G(E_T)$ on $D(X)$.}
\end{exe}
Notice that Proposition \ref{daniel} is not true if $E$ has no sink. Indeed, in  the next example, $\lambda$ is a group-type partial action of $\G(E)$ on $D(X)$ and $|d^{-1}(u)|=2,$ for some vertex $u$.

\begin{exe}\label{contra-exemplo} {\rm %In the previous examples, if $\lambda$ is a group-type partial action of $\G(E)$ on $D(X)$ and $|E^{0}|\geq 2$ then $|d^{-1}(v)|\leq 1$, for all $v\in E^{0}$. However, it is not true in general. In fact.
		Consider the directed graph $E$ given by 
		
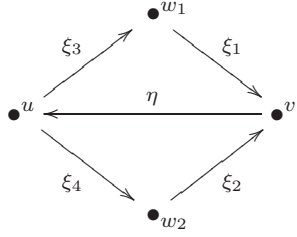
\begin{align*}
			\xymatrix{ &                                  &\hspace*{.4cm}\bullet^{w_1}\ar@{->}[rd]^{\xi_1}&  &  \\  
				&\bullet^{u}\ar@{->}[ru]^{\xi_3}\ar@{->}[rd]_{\xi_4} & &\bullet^{v} \ar@{->}[ll]_{\eta}& \\
				&                                                    &\hspace*{.4cm}\bullet_{w_2} \ar@{->}[ru]_{\xi_2} &                      &   }
		\end{align*}
		Consider the set $\tau(v)=\{\ta_v=v,\,\ta_{w_1}=\xi_1,\,\ta_u=\eta^{-1},\,\tau_{w_2}=\xi_2\}$. It is clear that $\ta_{w_1}$ and $\ta_{w_2}$ satisfy \eqref{eq1-tg} while $\ta_u$ verifies \eqref{eqq1-tg}. Hence, $\lambda$ is a group-type partial action of $\G(E)$ on $D(X)$ and $|d^{-1}(u)|=2$.}
\end{exe}

Recall that $E$ is {\it  finite} if both $E^0$ and $E^1$ are finite sets. 

\begin{remark}  \label{conse}{\rm If $E$  is finite and $E$ has no sinks, then $E$ contains a cycle. Thus $L_{\Bbbk}(E)$ is neither semisimple (Corollary  4.2.13 of \cite{AAM}) nor von Neumann regular  (Theorem 3.4.1 of \cite{AAM}).}
\end{remark}

\begin{exe}\label{LPA0}{\rm Let $n\in \N$, $n>1$ and consider the directed graph $E=A_n$:
		\begin{align*}
			\xymatrix{\bullet^{v_1}\ar@{->}[r]^{\xi_1}& \bullet^{v_2}\ar@{->}[r]^{\xi_2}& \bullet^{v_3}\ar@{.}[r]&\bullet^{v_{n-1}}\ar@{->}[r]^{\xi_{n-1}}&\bullet^{v}.} 	
		\end{align*}	
		%It is clear that $\G(E)$ is a connected groupoid. Take  $\ta(v)=\{\tau_{v_i}:v\to v_i\}$ a transversal for $v$ given by: 
		%\begin{align*}
		%&	\tau_{v}=v,& & \tau_{v_i}=\xi_i\xi_{i+1}\cdots \xi_{n-1},& &\text{for all }\,1\leq i \leq n-1.&
		%\end{align*}	
		%Notice that $\ta_{v_i}\in W^{(1)}$, for each $1\leq i\leq n-1$, and \eqref{eq1-tg} is satisfied. Thus,  $\lambda$ is a group-type partial action of $\G(E)$ on $D(X)$.  By  (b) in (ii) of Remark \ref{conse} %Observe that the isotropy group $\G(E)(v)=\{v\}$ is trivial. Hence,  Theorem \ref{teo-decomp} and Theorem \ref{iso-GY}  imply that
		%we have  $L_{\Bbbk}(E)\simeq D(X)\star_{\beta} (E^0)^2$, where $\beta=(B_u,\beta_u)_{u\in (E^0)^2}$ is the global action of $\G(E)_0^2=(E^0)^2$ on $D(X)$ given by \eqref{action-beta}. Explicitly, if $u=(v_i,v_j)\in (E^0)^2$ then $B_u=D(X)_{v_j}=D(X)1_{v_j}$. Since $X_{v_j}=\{\ta_{v_j}\}$, it follows that $1_{v_j}1_g=0$ if $\tau_{v_j}\notin X_g$ or $1_{v_j}1_g=1_{v_j}$ if $\tau_{v_j}\in X_g$, for all $g\in \G(E)$. Hence $D(X)1_{v_j}=\Bbbk\langle 1_{v_j}\rangle$ which is isomorphic (as algebra) to $\Bbbk$. Also, $\beta_u:\Bbbk\langle 1_{v_i}\rangle \to \Bbbk\langle 1_{v_j}\rangle$ is given by $\beta_u(a1_{v_i})=a1_{v_j}$, for all $a\in \Bbbk$. Hence it is clear that 
		%\begin{align*}
		%	\psi: D(X)\star_{\beta} (E^0)^2\to{\rm M}_n(\Bbbk),\quad \psi(a1_{v_j}\delta_{(v_i,v_j)})=ae_{ji},\quad a\in\Bbbk,
		%\end{align*}
		By Corollary \ref{cor-matrices}, we have that $L_{\Bbbk}(E)\simeq D(X)\star_{\beta} (E^0)^2\simeq{\rm M}_n(\Bbbk)$.}
\end{exe}

%{\color{blue} Inspired by Example \ref{LPA0} we give the next.
%\begin{prop1} Let $v\in E^0$ be a sink. Then $G(E)(v)$ is trivial, if and only if, for any $v_i\in E^0$ with $v_i\neq v$ there is a unique $\tau_{v_i}\in W^{(1)}$ such that $d(t_{v_i})=v_i$ and $r(v_i)=v.$ In this case  the transversal $\tau(v)=\{v, \tau_{v_i}\}_{i\in \Z^{+}}$ determines   group-type partial action of $\G(E)$ on $D(X)$ and $L_{\Bbbk}(E)\simeq D(X)\star_{\beta} E_0^2.$
%\end{prop1}
%\begin{proof} If  $G(E)(v)$ is not trivial there is $\xi\in G(E)(v)$ with $\xi\neq v,$  write $\xi=\xi_1\cdots \xi_n,$ for some $n>1, \xi_i\in E^1\cup (E^1)^*$ with $r(\xi_i)=d(\xi_{i+1}),$ then $d(\xi_1)=v=r(\xi_n),$ since $v$ is a sink $\xi_1\in (E^1)^*$  and $\xi_n\in E^1$ with $\xi_1\m \neq \xi_n,$ thus $\eta=\xi_1\m\in E^1$  \end{proof}}

\begin{exe}\label{LPA1}{\rm In the previous example, we take a transversal of a vertex $v$ which is a sink. Now we present an example where the vertex is not a sink. Consider the directed graph $E$:
		
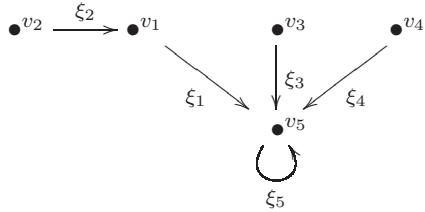
		\begin{align*}
			\xymatrix{ &\bullet^{v_2} \ar@{->}[r]^{\xi_2}&\bullet^{v_1}\ar@{->}[rd]_{\xi_1}&\bullet^{v_3}\hspace{-0.35cm}\ar@{->}[d]^{\xi_3}& \bullet^{v_4}\ar@{->}[ld]^{\xi_4} \\
				&                                 &                               &\hspace{0.35cm}\bullet^{v_5}\ar@(dl,dr)_{\xi_5}                & &} 
		\end{align*}
		Clearly, $\G(E)$ is connected. Let $\tau(v_1)=\{\tau_{v_i}:v_1\to v_i\}_{1\leq i\leq 5}$ be a transversal for $v_1$ given by
		\begin{align*}
			&\tau_{v_1}={v_1},& &\tau_{v_2}=\xi_2,& &\tau_{v_i}=\xi_{i}\xi\m_1,\,\,\, \text{for all } 3\leq i\leq 5.&	
		\end{align*}
		It is directly to verify that $\ta(v_1)\setminus\{v_1\}\subset W^{(1)}\cup W^{(2)}$ and \eqref{eq1-tg}, \eqref{eq2-tg} and \eqref{eq3-tg} are satisfied. Then, by Proposition \ref{gtype}, $\lambda$ is a group-type partial action of $\G(E)$ on $D(X)$. Notice that $\G(E)(v_1)$ is the infinite cyclic group generated by $\xi_1\xi_5\xi_1\m$ which is isomorphic (as a group) to $\Z$. By Theorem \ref{teo-decomp}, $L_{\Bbbk}(E)\simeq C\star_{\gamma}\Z$, where $C=D(X)\star_{\beta}(E^0)^2$.  We will calculate $C$ explicitly. Denote  by $\xi^{\infty}_5$ the infinite path $\xi_5\xi_5\xi_5\cdots $. Then $X_{v_i}=\{\xi_i\xi^{\infty}_5\}$ for all $i\neq 2$ and  $X_{v_2}=\{\xi_2\xi_1\xi^{\infty}_5\}$. Hence, as in Example \ref{LPA0}, we get $C\simeq {\rm M}_5(\Bbbk)  $ and $L_{\Bbbk}(E)\simeq {\rm M}_5(\Bbbk)\star_{\gamma}\Z$.}
\end{exe}

\begin{exe} {\rm Let $n$ be a positive integer and $R_n$ the \textit{rose with n petals},
		that is, $R_n$ is the directed graph that has one vertex $v$ and $n$ loops $\xi_i$, with $1\leq i\leq n$. Here $\G(E)=\G(E)(v)$ is a group. Thus the only transversal possible to $v$ is $\tau(v)=\{v\}$.  In this case,  the action $\beta$ of $\G(E)$ on $D(X)$ is trivial and $L_{\Bbbk}(R_n) \simeq D(X)\star_{\gamma}G(E)$. Moreover, $X=X_v=X_{\xi_i\m}=W^\infty$   and, by  Lemma \ref{decompo},  $D(X)=\bigoplus_{i=1}^nD(X)_{\xi_i}$.  It is straightforward to check that 
		\begin{align}\label{eq-ex}
			D(X)_{\xi\m_j}=D(X),\quad \lambda_{\xi_j}(1_{D(X)})=1_{\xi_j}\,\,\text{ and }\,\, \lambda_{\xi\m_j}(1_{\xi_i})=\delta_{i,j}1_{D(X)}.			
		\end{align}
		Since $G(E)_0=E^0=\{v\}$, it is easy to verify that the partial action $\gamma$ of $\G(E)_v$ on $D(X)$ given in Lemma 4.3 of \cite{BPP} is equal to $\lambda$. Thus $L_{\Bbbk}(R_n) \simeq D(X)\star_{\lambda}G(E)$. On the other hand, let $L_{\Bbbk}(1,n)$ be  the free associative $\Bbbk$-algebra with $2n$ generators $\{x_1,\ldots,x_n,y_1,\ldots,y_n\}$ subject to the relations $y_jx_i=\delta_{i,j}1_R$ and $\sum_{i=1}^{n}x_iy_i=1_R$. Consider the linear map $\varphi$ from the free associative algebra $\Bbbk\langle x_1,\ldots,x_n,y_1,\ldots,y_n\rangle$ onto $D(X)\star_{\lambda} \G(E)$ given by 
		\begin{align*}
			&\varphi(y_j)=1_{D(X)}\delta_{\xi\m_j}, & &\varphi(x_i)=1_{\xi_i}\delta_{\xi_i}.
		\end{align*}	
		By \eqref{eq-ex}, we have that
		\begin{align*}
			\varphi(y_j)\varphi(x_i)&=\big(1_{D(X)}\delta_{\xi\m_j}\big) \big(1_{\xi_i}\delta_{\xi_i}\big)=\delta_{i,j}1_{D(X)}\delta_{\xi\m_j\xi_i}=\delta_{i,j}1_{D(X)\star_{\lambda} \G(E)},\\
			\sum_{i=1}^{n}\varphi(x_i\varphi(y_i)&=\sum_{i=1}^{n}\big(1_{\xi_i}\delta_{\xi_i}\big)\big(1_{D(X)}\delta_{\xi\m_i}\big)=\sum_{i=1}^{n}1_{\xi_i}\delta_{1}=	1_{D(X)\star_{\lambda} \G(E)}.
		\end{align*}
		Hence, $\varphi$ induces an algebra homomorphism $\psi:L_{\Bbbk}(1,n)\to D(X)\star_{\lambda} \G(E)$ which is bijective. Then $L_{\Bbbk}(R_n) \simeq L_{\Bbbk}(1,n)$ and we recover Proposition 1.3.2 of \cite{AAM}. }
\end{exe}

\subsection{Ring theoretic properties of $L_{\Bbbk}(E)$}
Given a directed graph $E,$ we apply the results of Section 3   to relate algebraic properties of the ring $D(X)$ to  algebraic properties of the Leavitt path algebra $L_{\Bbbk}(E).$

\subsubsection{Noetherianity of $L_{\Bbbk}(E)$}  

We shall relate  noetherian properties  of  $L_{\Bbbk}(E)$ and $D(X).$
We start  with the following.

\begin{lem1}\label{noetd}   The ring $D(X)$ is  noetherian if and only if $E^0$ is finite and $D(X)_v$ is a  noetherian ring, for all $v\in E^0$.
	
\end{lem1}
\begin{proof}   
	Assume that $D(X)$ is  noetherian and suppose that  $E^0=\{v_0,v_1, \ldots\}$  is infinite. Consider, for each $j\geq 0$, the left ideal $I_j:=D(X)_{v_0}\oplus\ldots\oplus D(X)_{v_j}$ of $D(X)$. Since $1_{v_i}1_{v_j}=\delta_{ij}1_{v_i}$, we have that $1_{v_{j+1}}\in I_{j+1}\setminus I_j$. Hence, $I_j\subsetneq I_{j+1}$ for all $j\geq 0$. Thus we have a contradiction because $D(X)$ is noetherian. 
	Now consider $v\in E^0$ and $I$ a left ideal of $D(X)_v$. By Lemma \ref{decompo}, $I$ is a left ideal of $D(X)$ and whence $I$ is finitely generated as left $D(X)$-module. Then $I$ is finitely generated as left $D(X)_v$-module because  $D(X)_v$ is a direct summand of $D(X)$. The converse is immediate from well-known results of noetherian rings.
\end{proof}
{\it From now on we shall assume that  $\G(E)$ is connected,  $E^0$ is finite and that there are a vertex $v\in E^0$ and a transversal  $\tau(v)=\{\tau_{w}:v\to w\,:\, w\in E^0\}$ for $v$ such that $\ta(v)\setminus\{v\}\subset W^{(-1)}\cup W^{(1)}\cup W^{(2)}$ and \eqref{eqq1-tg}, \eqref{eq1-tg}, \eqref{eq2-tg} and \eqref{eq3-tg} are satisfied}.  \smallbreak

Observe that as  $\lambda$ is a group-type partial action of $\G(E)$ on $D(X)$ we obtain from Theorem \ref{teo-decomp} that  $L_{\Bbbk}(E)\simeq D(X)\star_{\lambda}\G(E)$. Also,  if $\G(E)(v)$ is finite then $\G(E)$ is finite by Proposition \ref{group:connec}.

\smallbreak

For the reader's convenience, we recall some extra notions on a directed graph $E$.  Let  $\xi=\xi_1\xi_2\ldots \xi_n$ be a path in $E$. An element $\eta\in E^1$ is called an {\it exit for $\xi$} if there exists $1\leq i\leq n$ such that $d(\eta)=d(\xi_i)$ and  $\eta\neq \xi_l$ and  we say that $E$ satisfies {\it Condition (NE)} if there is no cycle in $E$ with an exit.

\begin{prop1}\label{noethercond}  The following assertions are valid.
	\begin{enumerate}[\rm (i)]
		\item If  $L_{\Bbbk}(E)$ is left noetherian then $D(X)$ and $D(X)\star_{\bt} (E^0)^2$ are left noetherian.\smallbreak
		
		\item Suppose that $G(E)(v)$ is a polycyclic-by-finite group. The following statements are equivalent:
		\begin{enumerate} [\rm (a)]
			\item $L_{\Bbbk}(E)$ is left noetherian,\smallbreak
			\item $D(X)$  is left noetherian,\smallbreak
			\item  $D(X)_v$ is left noetherian, for all $v\in E^0$,\smallbreak
			\item  $E$ is finite and it satisfies Condition  (NE).
		\end{enumerate}
	\end{enumerate}
\end{prop1}

\begin{proof} 
	Item (i) follows directly from Theorem  \ref{noete} (i).
	The equivalence (a) $\Leftrightarrow$ (d) in 	(ii)  follows from Theorem 3.10 of \cite{AAM2} and  (b) $\Leftrightarrow$ (c) is Lemma \ref{noetd}.  From Theorem \ref{noete} (ii) we obtain that (b) $\Rightarrow$ (a). Notice that  (a) $\Rightarrow$ (b) by item (i).
\end{proof}

\begin{remark}\label{obs-true-general} {\rm
		The equivalence between (a) and (d) in (ii) of Proposition \ref{noethercond}  was proved for general directed graphs in Theorem 3.10 of \cite{AAM2}.	}
\end{remark}

\begin{remark}\label{rem-advantage} {\rm Since $\Z$ is a polycyclic-by-finite and  ${\rm M}_5(\Bbbk)$ is left noetherian it follows from Proposition \ref{noethercond} (ii) that  the algebra $L_{\Bbbk}(E)$ from Example \ref{LPA1}  is left noetherian. If we realize $L_{\Bbbk}(E)$ as a partial skew group ring as in Proposition 3.2 of  \cite{GR}, we obtain that $L_{\Bbbk}(E)$ is isomorphic to the partial skew group ring $D_0\star_{\theta}\mathbb{F}_5$, where $\mathbb{F}_5$ is the free group generated by $E^{1}$, $D_0$ is a subalgebra of the algebra of functions on a set $X$ and $\theta$ is a partial action of $\mathbb{F}_5$ on $D_0$. It was proved in \cite{CCF} that a partial skew group ring of a noetherian ring by a partial action of a polycyclic-by-finite group is noetherian. It is well-known that $\mathbb{F}_5$ is not polycyclic-by-finite. Hence we can not use \cite{CCF} for to conclude that $L_{\Bbbk}(E)$ is noetherian. }
\end{remark}

\subsubsection{Regularity of $L_{\Bbbk}(E)$}

There are several results in the literature on group  graded von Neumann regularity for Leavitt path algebras; see, for instance, Theorems 1.1, 1.3 and 1.4 in \cite{La2}. In particular,  $L_{\Bbbk}(E)$ with the standard $\Z$-grading is graded von Neumann regular.  Our next result relates groupoid graded von Neumann regularity of $L_{\Bbbk}(E)$ with von Neumann regularity of the rings $D(X)$ and $L_{\Bbbk}(E)$ itself.

\begin{prop1} 
	Suppose that $\operatorname{ch}(\Bbbk)$ does not divide  $\vert E^0 \vert1_{\Bbbk}$. 
	\begin{enumerate}[\rm (i)]
		\item The following statements are equivalent:
		\begin{enumerate}[\rm (a)]
			\item $D(X)$ is von Neumann regular,\smallbreak
			\item $D(X)\star_{\bt}\G_0^2$ is graded  von Neumann regular,\smallbreak
			\item $D(X)\star_{\bt}\G_0^2$ is von Neumann regular,\smallbreak
			\item $L_{\Bbbk}(E)$ is  $\G(E)$-graded von Neumann regular.\smallbreak
		\end{enumerate}

		\item If $\G(E)(v)$ is finite and $ {\rm tr}_{\lambda_{\G(E)(v)}}(1_v)$ is invertible in $D(X)_v$ then:\smallbreak
		\begin{enumerate}[\rm (a)]
			\item  The ring $D(X)$ is semiprimitive, that is $J(D(X))=\{0\}$.\smallbreak
			\item If $L_{\Bbbk}(E)$ is  $\G(E)$-graded von Neumann regular then it is
			von Neumann regular,\smallbreak
		\end{enumerate}
	\end{enumerate}
	\begin{proof} Firstly, notice that $\vert E^0 \vert 1_{\Bbbk}$ is invertible in $D(X)$.
		Also, by Theorem \ref{teo-decomp}, we have that $L_{\Bbbk}(E)\simeq D(X)*_{\lambda} \G(E)$ because $\lambda$ is a group-type partial action of $\G(E)$ on $D(X)$. Thus (i) follows from   Theorem \ref{vn1}.\smallbreak

		\noindent (ii) From Theorem \ref{jaco}, $J(L_{\Bbbk}(E))\simeq J(D(X))*_\lambda \G(E)$. However, by Proposition 6.3 of \cite{AAP}, $J(L_{\Bbbk}(E))=\{0\}$ and whence $J(D(X))=\{0\}$. Hence $D(X)$ is semiprimitive. For (b), observe that from the equivalence (a) $\Leftrightarrow$ (d) in  part (i) we obtain that $D(X)$ is von Neumann regular.  Hence (b) follows by  (e) in Theorem  \ref{vn1}.
	\end{proof}
\end{prop1}

\begin{exe}{\rm  Assume that $\operatorname{ch}(\Bbbk)\neq 5$ and consider the Leavitt path algebra given in Example \ref{LPA1}, that is, $L_{\Bbbk}(E)\simeq {\rm M}_5(\Bbbk)\star_{\gamma}\Z$. Since ${\rm M}_5(\Bbbk)$ is von Neummann regular, it follows from Proposition 3.10 of \cite{La2} and Remark \ref{obs-g-grade} that $L_{\Bbbk}(E)$ is $\G(E)$-graded von Neumann regular. Hence, the equivalences in (i) of the previous proposition are satisfied. Also, Remark \ref{obs-g-grade} implies that $L_{\Bbbk}(E)$ is $\G(x)$-graded von Neumann. However, by Theorem 3.4.1 of \cite{AAM}, $L_{\Bbbk}(E)$ is not von Neumann regular because $E$ has one cycle.}
\end{exe}

\begin{remark}{\rm We have  by (i) of Lemma \ref{decompo} that $D(X)$ is a direct sum of orthogonal idempotents rings. Thus, by Lemma \ref{vnr},  $D(X)$ is von Neumann regular if and only if $D(X)_v$ is  von Neumann regular, for all $v\in E^0$. }
\end{remark}

\subsubsection{Right perfect, semiprimary and Krull dimension}

We end the paper with the following remark about other ring theoretic properties of $L_{\Bbbk}(E)$.

\begin{remark}{\rm
		\begin{enumerate}[\rm (i)]
			\item  $L_{\Bbbk}(E)$ is right perfect if and only if $D(X)$ is right perfect and $\G(E)(v)$ is finite; this follows from Theorem \ref{rightp}.\smallbreak
			
			\item  $L_{\Bbbk}(E)$ is semiprimary if and only if $D(X)$ is semiprimary and $\G(E)(v)$ is finite; this follows from Theorem \ref{spp}\smallbreak
			
			\item  	If $\G(E)(v)$ is finite and $D(X)$ has right Krull dimension then $L_{\Bbbk}(E)$ has right Krull dimension and $\operatorname{Kdim} L_{\Bbbk}(E)\leq \operatorname{Kdim} D(X)$.
			In fact, this is an immediate consequence of Theorem \ref{krull} since $\G(E)$ is finite. 
	\end{enumerate}}
\end{remark}

%\begin{acknowledgements}
%If you'd like to thank anyone, place your comments here
%and remove the percent signs.
%\end{acknowledgements}

% Authors must disclose all relationships or interests that 
% could have direct or potential influence or impart bias on 
% the work: 
%
% \section*{Conflict of interest}
%
% The authors declare that they have no conflict of interest.

% BibTeX users please use one of
%\bibliographystyle{spbasic}      % basic style, author-year citations
%\bibliographystyle{spmpsci}      % mathematics and physical sciences
%\bibliographystyle{spphys}       % APS-like style for physics
%\bibliography{}   % name your BibTeX data base

% Non-BibTeX users please use

\end{document}